\RequirePackage{fix-cm}
\documentclass[smallextended]{svjour3}       
\smartqed  
\usepackage{graphicx}
\usepackage{amsmath,amssymb}
\usepackage{array}
\usepackage{enumitem}
\usepackage{xcolor}
\usepackage{bm}
\usepackage[normalem]{ulem}

\providecommand{\etal}		{\emph{et al\@.}\xspace}
\providecommand{\ie}		{\emph{i.e\@.}\xspace}

\providecommand{\myurl}[1][]	{\texttt{web.eecs.umich.edu/$\sim$fessler#1}\xspace}

\providecommand{\onweb}[1]	{Available from \myurl.}

\long\def\comment#1{}

\providecommand{\bcent}		{\begin{center}}
\providecommand{\ecent}		{\end{center}}
\providecommand{\benum}		{\begin{enumerate}}
\providecommand{\eenum}		{\end{enumerate}}
\providecommand{\bitem}		{\begin{itemize}}
\providecommand{\eitem}		{\end{itemize}}
\providecommand{\bvers}		{\begin{verse}}
\providecommand{\evers}		{\end{verse}}
\providecommand{\btab}		{\begin{tabbing}}	
\providecommand{\etab}		{\end{tabbing}}

\newcounter{blist}
\providecommand{\blistmark}	{\makebox[0pt]{$\bullet$}}
\providecommand{\blistitemsep}	{0pt}
\providecommand{\blist}[1][]	{%
\begin{list}{\blistmark}{%
\usecounter{blist}%
\setlength{\itemsep}{\blistitemsep}%
\setlength{\parsep}{0pt}%
\setlength{\parskip}{0pt}%
\setlength{\partopsep}{0pt}%
\setlength{\topsep}{0pt}%
\setlength{\leftmargin}{1.2em}%
\setlength{\labelsep}{0.5\leftmargin}
\setlength{\labelwidth}{0em}%
#1}
}
\providecommand{\elist}		{\end{list}}

\providecommand{\blistitemsep}	{0pt}
\providecommand{\bjfenum}[1][]	{%
\begin{list}{\bcolor{\arabic{blist}.} }{%
\usecounter{blist}%
\setlength{\itemsep}{\blistitemsep}%
\setlength{\parsep}{0pt}%
\setlength{\parskip}{0pt}%
\setlength{\partopsep}{0pt}%
\setlength{\topsep}{0pt}%
\setlength{\leftmargin}{0.0em}%
\setlength{\labelsep}{1.0\leftmargin}
\setlength{\labelwidth}{0pt}%
#1}
}

\newcounter{blistAlph}
\providecommand{\blistAlph}[1][]
{\begin{list}{\makebox[0pt][l]{\Alph{blistAlph}.}}{%
\usecounter{blistAlph}%
\setlength{\itemsep}{0pt}\setlength{\parsep}{0pt}%
\setlength{\parskip}{0pt}\setlength{\partopsep}{0pt}%
\setlength{\topsep}{0pt}%
\setlength{\leftmargin}{1.2em}%
\setlength{\labelsep}{1.0\leftmargin}
\setlength{\labelwidth}{0.0\leftmargin}#1}%
}

\newcounter{blistRoman}
\providecommand{\blistRoman}[1][]
{\begin{list}{\Roman{blistRoman}.}{%
\usecounter{blistRoman}%
\setlength{\itemsep}{0.5em}\setlength{\parsep}{0pt}%
\setlength{\parskip}{0pt}\setlength{\partopsep}{0pt}%
\setlength{\topsep}{0pt}%
\setlength{\leftmargin}{4em}%
\setlength{\labelsep}{0.4\leftmargin}
\setlength{\labelwidth}{0.6\leftmargin}#1}%
}

\usepackage{bbm} 
\providecommand{\qed}[1][0pt]	{\hfill\raisebox{#1}{\inmath{\Box}}} 

\providecommand{\defequ}	{\stackrel{\bigtriangleup}{=}}

\providecommand{\floor}[1]	{\xmath{\left\lfloor #1 \right\rfloor}}
\providecommand{\ceil}[1]	{\xmath{\left\lceil #1 \right\rceil}}

\let\equivsave\equiv
\def\equiv{\xmath{\equivsave}}

\providecommand{\ba}[1]		{\left[ \begin{array}{#1}}
\providecommand{\ea}		{\end{array} \right]}
\providecommand{\be}		{\begin{equation}}
\providecommand{\ee}[1]		{\label{#1}\end{equation}}
\providecommand{\bea}		{\begin{eqnarray}}
\providecommand{\eea}[1]	{\label{#1}\end{eqnarray}}
\providecommand{\beas}		{\begin{eqnarray*}}
\providecommand{\eeas}		{\end{eqnarray*}}
\providecommand{\beals}[1][1]	{\begin{alignat*}{#1}}	
\providecommand{\eeals}		{\end{alignat*}}

\providecommand{\berr}[2]{
\bgroup
\renewcommand{\theequation}{#1}
\be
#2
\ee{e,#1}
\egroup
\ignorespaces
}

\providecommand{\bearr}[2]{
\bgroup
\renewcommand{\theequation}{#1}
\bea
#2
\eea{e,#1}
\egroup
\ignorespaces
}

\providecommand{\inmath}	{\ensuremath}
\providecommand{\xmath}[1]	{\inmath{#1}\xspace}
\providecommand{\bmath}[1]	{\xmath{\bm{#1}}}	

\providecommand{\paren}[1]	{\xmath{\left(#1\right)}}

\providecommand{\braces}[1]	{\xmath{\left\{#1\right\}}}

\providecommand{\xfrac}[2]	{\xmath{\frac{#1}{#2}}}

\providecommand{\Frac}[2]	{\xmath{{#1}/{#2}}}
\providecommand{\pfrac}[3][]	{\xmath{\!\paren{#1\frac{#2}{#3}}}}

\providecommand{\wordbrace}[3][]{\,\inmath{\mathrm{#2}#1\!\braces{#3}}}

\providecommand{\diag}[1]	{\wordbrace{diag}{#1}}

\providecommand{\sinc}		{\Ofun{\mathrm{sinc}}}

\usepackage{jf-accent}
\usepackage{jf-fun}
\usepackage{jf-names}

\newcommand{\st} {\xmath{\text{s.t.}\:}}

\renewcommand{\defequ} {\triangleq}

\renewcommand{\defequ} {\triangleq}

\newcommand{\DT} {DT}
\newcommand{\primary} {secondary}
\newcommand{\auxiliary} {primary}
\newcommand{\dd} {\xmath{d}}

\newcommand{\bF} {\xmath{\breve{f}}}

\newcommand{\tS} {\bmath{Q}} 
\newcommand{\ts} {\bmath{q}} 
\newcommand{\tw} {\bmath{{\check{w}}}}
\newcommand{\ttau} {\bmath{{\check{\tau}}}}

\newcommand{\Sr} {\bmath{\breve{S}}}
\newcommand{\Qr} {\bmath{\breve{Q}}}
\newcommand{\qr} {\bmath{\breve{q}}}

\newcommand{\bmw} {\bmath{w}}
\newcommand{\bmlam} {\bmath{\lambda}}
\newcommand{\bmtau} {\bmath{{\tau}}}
\newcommand{\bmr} {\bmath{r}}
\newcommand{\bmh} {\bmath{h}}
\newcommand{\bmg} {\bmath{g}}
\newcommand{\bmG} {\bmath{G}}
\newcommand{\bmdel} {\bmath{\delta}}
\newcommand{\bmu} {\bmath{u}}
\newcommand{\tbmu} {\bmath{\check{u}}}

\newcommand{\bme} {\bmath{e}}
\newcommand{\bmx} {\bmath{x}}
\newcommand{\bmy} {\bmath{y}}
\newcommand{\bmz} {\bmath{z}}
\newcommand{\bmnu} {\bmath{\nu}}

\newcommand{\bmS} {\bmath{S}}

\newcommand{\bmA} {\bmath{A}}

\newcommand{\bmD} {\bmath{D}}

\newcommand{\bmt} {\bmath{t}}
\newcommand{\bmtt} {\bmath{\check{t}}}
\newcommand{\bmttheta} {\bmath{\check{\theta}}}
\newcommand{\bh} {\bmath{\bar{h}}}
\newcommand{\blam} {\xmath{\bar{\lambda}}}
\newcommand{\btau} {\xmath{\bar{\tau}}}
\newcommand{\bgam} {\xmath{\bar{\gamma}}}

\newcommand{\bmblam} {\bmath{\bar{\lambda}}}
\newcommand{\bmbtau} {\bmath{{\bar{\tau}}}}

\newcommand{\hlam} {\xmath{\hat{\lambda}}}
\newcommand{\htau} {\xmath{\hat{\tau}}}
\newcommand{\hgam} {\xmath{\hat{\gamma}}}

\newcommand{\bmhlam} {\bmath{\hat{\lambda}}}
\newcommand{\bmhtau} {\bmath{{\hat{\tau}}}}

\newcommand{\cF} {\xmath{\mathcal{F}_L(\Reals^d)}}
\newcommand{\cL} {\xmath{\mathcal{L}}}

\newcommand{\Reals} {\xmath{\mathbb{R}}}

\newcommand{\rki} {\xmath{r_{i,k}}}
\newcommand{\rkN} {\xmath{r_{N,k}}}

\newcommand{\rNmN} {\xmath{r_{N,N-1}}}

\newcommand{\bhzo} {\xmath{\bar{h}_{1,0}}}

\newcommand{\hki} {\xmath{h_{i,k}}}
\newcommand{\hkip} {\xmath{h_{i+1,k}}}

\newcommand{\hkj} {\xmath{h_{j,k}}}

\newcommand{\bhkn} {\xmath{\bar{h}_{n,k}}}
\newcommand{\bhknp} {\xmath{\bar{h}_{n+1,k}}}

\newcommand{\bhnnp} {\xmath{\bar{h}_{n+1,n}}}

\newcommand{\Hh} {\bmath{\hat{h}}}
\newcommand{\Hr} {\bmath{\hat{r}}}
\newcommand{\Hhzo} {\xmath{\hat{h}_{1,0}}}

\newcommand{\Hrki} {\xmath{\hat{r}_{i,k}}}
\newcommand{\Hhki} {\xmath{\hat{h}_{i,k}}}
\newcommand{\Hhkip} {\xmath{\hat{h}_{i+1,k}}}
\newcommand{\Hhkj} {\xmath{\hat{h}_{j,k}}}
\newcommand{\Hhkn} {\xmath{\hat{h}_{n,k}}}
\newcommand{\Hhknp} {\xmath{\hat{h}_{n+1,k}}}

\newcommand{\Hhiip} {\xmath{\hat{h}_{i+1,i}}}
\newcommand{\Hhimi} {\xmath{\hat{h}_{i,i-1}}}
\newcommand{\Hhimip} {\xmath{\hat{h}_{i+1,i-1}}}
\newcommand{\Hhnmn} {\xmath{\hat{h}_{n,n-1}}}
\newcommand{\Hhnnp} {\xmath{\hat{h}_{n+1,n}}}
\newcommand{\Hhnmnp} {\xmath{\hat{h}_{n+1,n-1}}}

\newcommand{\bhki} {\xmath{\bar{h}_{i,k}}}
\newcommand{\bhkip} {\xmath{\bar{h}_{i+1,k}}}
\newcommand{\bhkj} {\xmath{\bar{h}_{j,k}}}
\newcommand{\bhiip} {\xmath{\bar{h}_{i+1,i}}}
\newcommand{\bhimi} {\xmath{\bar{h}_{i,i-1}}}
\newcommand{\bhimip} {\xmath{\bar{h}_{i+1,i-1}}}

\numberwithin{equation}{section}

\usepackage[b5paper]{geometry}
\usepackage{fullpage}
\usepackage{hyperref}
\hypersetup{
        linktoc=page,
        colorlinks=true,
        linkcolor=blue,
        citecolor=red,
        filecolor=magenta,
        urlcolor=cyan
}

\begin{document}

\title{Optimized first-order methods for smooth convex minimization
\thanks{This research was supported in part 
by NIH grant R01-HL-098686
and U01-EB-018753.}
\thanks{After this paper was finalized for publication at 
\url{http://dx.doi.org/10.1007/s10107-015-0949-3}, 
we found a small error.  
This is a corrected version with deletions marked in {\color{red} red} 
and additions in {\color{blue} blue}.
}
}


\author{Donghwan Kim         \and
        Jeffrey A. Fessler 
}


\institute{Donghwan Kim \and Jeffrey A. Fessler \at
              Dept. of Electrical Engineering and Computer Science,
		University of Michigan, Ann Arbor, MI 48109 USA \\
              \email{kimdongh@umich.edu, fessler@umich.edu}           
}

\date{Date of current version: \today} 

\maketitle

\begin{abstract}
We introduce new optimized first-order methods
for smooth unconstrained convex minimization.
Drori and Teboulle~\cite{drori:14:pof} recently described a numerical method 
for computing the $N$-iteration optimal step coefficients 
in a class of first-order algorithms
that includes gradient methods,
heavy-ball methods~\cite{polyak:64:smo},
and Nesterov's fast gradient methods
\cite{nesterov:83:amo,nesterov:05:smo}.
However, the numerical method in~\cite{drori:14:pof}
is computationally expensive for large $N$,
and the corresponding numerically optimized 
first-order algorithm in~\cite{drori:14:pof}
requires impractical memory and computation
for large-scale optimization problems.
In this paper,
we propose optimized first-order algorithms
that achieve a convergence bound
that is two times smaller than for Nesterov's fast gradient methods;
our bound is found analytically and refines the numerical bound in~\cite{drori:14:pof}.
Furthermore, the proposed optimized first-order methods
have efficient forms
that are remarkably similar to Nesterov's fast gradient methods.

\end{abstract}

\section{Introduction}
\label{intro}

First-order algorithms are used widely to solve large-scale optimization problems
in various fields such as signal and image processing, 
machine learning, communications and many other areas.
The computational cost per iteration of first-order algorithms
is mildly dependent on the dimension of the problem,
yielding computational efficiency.
Particularly, Nesterov's fast gradient methods
\cite{nesterov:83:amo,nesterov:05:smo}
have been celebrated in various applications
for their fast convergence rates and efficient implementation.
This paper proposes first-order algorithms
(OGM1 and OGM2 in Section~\ref{sec:opt,fo,form})
that achieve a worst-case convergence bound that is twice as small
as Nesterov's fast gradient methods
for smooth unconstrained convex minimization
yet have remarkably similar 
efficient implementations.

We consider finding a minimizer over $\Reals^d$
of a cost function $f$ belonging to
a set \cF of smooth convex functions
with $L$-Lipschitz continuous gradient.
The class of first-order (FO) algorithms of interest generates a sequence of 
points $\{\bmx_i\in\Reals^d\;:\;i=0,\cdots,N\}$ 
using the following scheme:

\fbox{
\begin{minipage}[t]{0.85\linewidth}
\vspace{-10pt}
\begin{flalign}
&\quad \text{\bf Algorithm Class FO} & \nonumber \\
&\qquad \text{Input: } f\in \cF,\; \bmx_0\in\Reals^d. & \nonumber \\
&\qquad \text{For } i = 0,\cdots,N-1 & \nonumber \\
&\qquad \qquad \bmx_{i+1} = \bmx_i - \frac{1}{L}\sum_{k=0}^i \hkip f'(\bmx_k). & \label{eq:fo}
\end{flalign}
\end{minipage}
} \vspace{5pt}

\noindent
The update step at the $i$th iterate $\bmx_i$ uses
a linear combination of previous and current gradients
$\{f'(\bmx_0),\cdots,f'(\bmx_i)\}$.
The coefficients $\{\hki\}_{0\le k<i\le N}$ determine the step size
and are selected prior to iterating (non-adaptive).
Designing these coefficients appropriately is the key
to establishing fast convergence.
The algorithm class FO includes gradient methods,
heavy-ball methods~\cite{polyak:64:smo},
Nesterov's fast gradient methods~\cite{nesterov:83:amo,nesterov:05:smo},
and our proposed optimized first-order methods.

Evaluating the convergence bound of 
such first-order algorithms is essential.
Recently, Drori and Teboulle (hereafter ``\DT'')~\cite{drori:14:pof} considered
the Performance Estimation Problem (PEP) approach to 
bounding the decrease of a cost function $f$.
For given coefficients $\bmh = \{\hki\}_{0\le k<i\le N}$,
a given number of iterations $N\ge1$
and a given upper bound $R>0$ on
the distance between an initial point $\bmx_0$ 
and an optimal point $\bmx_* \in X_*(f)\defequ \argmin{\bmx\in\Reals^d}f(\bmx)$,
the worst-case performance bound of a first-order method 
over all smooth convex functions $f\in\cF$
is the solution of
the following constrained optimization problem\footnote{
\label{ft}
The problem $\mathcal{B}_{\mathrm{P}}(\bmh,N,\dd,L,R)$
was shown to be independent of $d$ 
{\color{blue}
for the large-scale condition ``$d \ge N+2$''
}
in~\cite{taylor:15:ssc};
{\color{red} \sout{thus}}
this paper's results are independent of $d$
{\color{blue}
for any $d \ge 1$}
.
}
\cite{drori:14:pof}:
\begin{align}
\mathcal{B}_{\mathrm{P}}(\bmh,N,\dd,L,R) \defequ\;
&\max_{f\in\cF}  
	\max_{\substack{\bmx_0,\cdots,\bmx_N\in\Reals^d, \\ 
		\bmx_*\in X_*(f)}}
	f(\bmx_N) - f(\bmx_*) 
\label{eq:PEP} \tag{P} \\
&\st \; \bmx_{i+1} = \bmx_i - \frac{1}{L}\sum_{k=0}^i \hkip f'(\bmx_k), 
	\quad i=0,\cdots,N-1, \nonumber \\
&\quad\;\;\;
	||\bmx_0 - \bmx_*|| \le R. \nonumber 
\end{align}
As reviewed in Section~\ref{sec:pep},
\DT~\cite{drori:14:pof} used relaxations to
simplify the intractable problem~\eqref{eq:PEP}
to a solvable form.

Nesterov's fast gradient methods
\cite{nesterov:83:amo,nesterov:05:smo}
achieve the optimal rate of decrease $O\paren{\frac{1}{N^2}}$
for minimizing a smooth convex function $f$
\cite{nesterov:04}.
Seeking first-order algorithms that converge faster
(in terms of the constant factor)
than Nesterov's fast gradient methods,
\DT~\cite{drori:14:pof}
proposed using a (relaxed) PEP approach 
to optimize the choice of \bmh in class FO
by minimizing a (relaxed) worst-case bound at the $N$th iteration
with respect to \bmh.
In~\cite{drori:14:pof}, the optimized \bmh factors were computed numerically,
and were found to yield faster convergence than Nesterov's methods.
However, numerical optimization of \bmh in~\cite{drori:14:pof}
becomes expensive for large $N$.
In addition, the general class FO 
requires $O(N^2d)$ arithmetic operations 
for $N$ iterations
and $O(Nd)$ memory for storing all gradients 
$\{f'(\bmx_i) \in \Reals^d\;:\;i=0,\cdots,N-1\}$,
which is impractical for large-scale problems.

This paper proposes optimized first-order algorithms that
have a worst-case convergence bound
that is twice as small as that of
Nesterov's fast gradient methods,
inspired by \DT~\cite{drori:14:pof}.
We develop remarkably efficient formulations
of the optimized first-order algorithms
that resemble those of Nesterov's fast gradient methods,
requiring $O(Nd)$ arithmetic operations  
and $O(d)$ memory.

Section~\ref{sec:prob} reviews the smooth convex minimization problem
and introduces the approach to optimizing \bmh
used here and in~\cite{drori:14:pof}.
Section~\ref{sec:fo,ex} illustrates 
Nesterov's fast gradient methods that are in class FO.
Section~\ref{sec:pep} reviews \DT's (relaxed) PEP approach
and Section~\ref{sec:fgm,pep,ana} uses it to derive a new convergence bound 
for the \emph{\primary}~variables in Nesterov's fast gradient methods.
Section~\ref{sec:opt,fo} reviews \DT's analysis
on numerically optimizing \bmh using (relaxed) PEP
for first-order methods,
and derives an analytical form of the optimized coefficients \bmh
and a corresponding new analytical bound.
Section~\ref{sec:opt,fo,form}
investigates efficient formulations of the proposed first-order methods
(OGM1 and OGM2).
Section~\ref{sec:disc} shows that 
the corresponding analytical upper bound is tight
and Section~\ref{sec:conc} concludes.

\section{Problem and approach}
\label{sec:prob}

\subsection{Smooth convex minimization problem}

We consider first-order algorithms 
for solving the following minimization problem
\begin{align}
\min_{\bmx\in\Reals^d} \;&\; f(\bmx)
\label{eq:prob} \tag{M}
,\end{align}
where the following two conditions are assumed:
\begin{itemize}[leftmargin=40pt]
\item  
$f\;:\;\Reals^d\rightarrow\Reals$ is 
a convex function of the type $C_L^{1,1}(\Reals^d)$, 
\ie, continuously differentiable with Lipschitz continuous gradient: 
\begin{align*}
||f'(\bmx) - f'(\bmy)|| \le L||\bmx-\bmy||, \quad \forall \bmx, \bmy\in\Reals^d
,\end{align*}
where $L > 0$ is the Lipschitz constant.
\item 
The optimal set $X_*(f)=\argmin{\bmx\in\Reals^d} f(\bmx)$ is nonempty,
\ie, the problem~\eqref{eq:prob} is solvable.
\end{itemize}

We focus on measuring the ``inaccuracy'' $f(\bmx_N) - f(\bmx_*)$ after $N$ iterations
to quantify the worst-case performance of any given first-order algorithm.

\subsection{Optimizing the step coefficients \bmh
of first-order algorithms}

In search of the best-performing first-order methods,
\DT~\cite{drori:14:pof} 
proposed to optimize $\bmh = \{\hki\}_{0\le k < i \le N}$ in Algorithm FO 
by minimizing a worst-case bound of $f(\bmx_N) - f(\bmx_*)$ 
for a given number of iterations $N\ge1$
and initial distance $R > 0$,
by adding $\argmin{\bmh}$ to problem~\eqref{eq:PEP} as follows:
\begin{align}
\Hh_{\mathrm{P}} \defequ
\argmin{\bmh\in\Reals^{\Frac{N(N+1)}{2}}}
\mathcal{B}_{\mathrm{P}}(\bmh,N,\dd,L,R)
\label{eq:HP} \tag{HP} 
.\end{align}
Note that $\Hh_{\mathrm{P}}$ is independent\footnote{ 
Substituting $\bmx' = \frac{1}{R}\bmx$ 
and $\bF(\bmx') = \frac{1}{LR^2}f(R\bmx')\in\mathcal{F}_1(\Reals^d)$ 
in problem~\eqref{eq:PEP}, 
we get
$\mathcal{B}_{\mathrm{P}}(\bmh,N,L,R) = LR^2\mathcal{B}_{\mathrm{P}}(\bmh,N,1,1)$.
This leads to
$\Hh_{\mathrm{P}} = \argmin{\bmh} \mathcal{B}_{\mathrm{P}}(\bmh,N,L,R)
        = \argmin{\bmh} \mathcal{B}_{\mathrm{P}}(\bmh,N,1,1)$.
}
of 
$L$ and $R$.
{\color{red} \sout{and $d$. (See footnote~\ref{ft}.)}}
Solving problem~\eqref{eq:HP} would give
the step coefficients of the optimal first-order algorithm
achieving the best worst-case convergence bound.
\DT~\cite{drori:14:pof} relaxed\footnote{
Using the term `best' or `optimal' here
for~\cite{drori:14:pof} may be too strong,
since~\cite{drori:14:pof} 
relaxed~\eqref{eq:HP} to a solvable form.
We also use these relaxations,
so we use the term ``optimized'' for our proposed algorithms.
}
problem~\eqref{eq:HP}
to a tractable form,
as reviewed in Sections~\ref{sec:pep} and~\ref{sec:opt,fo,num}.
After these simplifications,
the resulting solution was computed
by a semidefinite program (SDP)
that remains computationally expensive 
for large $N$~\cite{drori:14:pof}.
In addition, the corresponding numerically optimized first-order algorithm
was impractical for large-scale problems,
requiring a linear combination of previous and current gradients
$\{f'(\bmx_0),\cdots,f'(\bmx_i)\}$
at the $(i+1)$-th iteration.\footnote{
If coefficients \bmh in Algorithm FO have a special recursive form,
it is possible to find an equivalent efficient form,
as discussed in Sections~\ref{sec:fo,ex} and~\ref{sec:opt,fo,form}.
}

To make \DT's work~\cite{drori:14:pof} practical,
we directly derive the ``analytical'' solution for \bmh
in a relaxed version of the problem~\eqref{eq:HP},
circumventing the numerical approach in~\cite{drori:14:pof}.
Interestingly, the analytical solution of 
the relaxed version of~\eqref{eq:HP}
satisfies a convenient recursion,
so we provide practical optimized algorithms
similar to Nesterov's efficient fast gradient methods.

\section{Nesterov's fast gradient methods}
\label{sec:fo,ex}

This section reviews 
Nesterov's well-known fast gradient methods~\cite{nesterov:83:amo,nesterov:05:smo}.
We further show the equivalence\footnote{ 
The equivalence of two of Nesterov's fast gradient methods
for smooth unconstrained convex minimization
was previously mentioned without details in~\cite{tseng:10:aag}.
}
of two of Nesterov's fast gradient methods
in smooth unconstrained convex minimization.
The analysis techniques used here 
will be important in Section~\ref{sec:opt,fo,form}.

\subsection{Nesterov's fast gradient method 1 (FGM1)}

Nesterov's first fast gradient method is
called FGM1~\cite{nesterov:83:amo}:

\fbox{
\begin{minipage}[t]{0.85\linewidth}
\vspace{-10pt}
\begin{flalign}
&\quad \text{\bf Algorithm FGM1} & \nonumber \\
&\qquad \text{Input: } f\in C_L^{1,1}(\Reals^d)\text{ convex},\; \bmx_0\in\Reals^d,\;
	\bmy_0 = \bmx_0,\; t_0 = 1. & \nonumber \\
&\qquad \text{For } i = 0,\cdots,N-1 & \nonumber \\
&\qquad \qquad \bmy_{i+1} = \bmx_i - \frac{1}{L}f'(\bmx_i) & \nonumber \\
&\qquad \qquad t_{i+1} = \frac{1+\sqrt{1+4t_i^2}}{2} & \nonumber \\
&\qquad \qquad \bmx_{i+1} 
	= \bmy_{i+1} + \frac{t_i-1}{t_{i+1}}(\bmy_{i+1} - \bmy_i). & \label{eq:ti}
\end{flalign}
\end{minipage}
} \vspace{5pt}

\noindent
Note that $t_i$ in~\eqref{eq:ti} satisfies 
the following relationships 
used frequently in later derivations:
\begin{align}
t_{i+1}^2 - t_{i+1} - t_i^2 = 0, \quad
t_i^2 = \sum_{k=0}^i t_k, \quad \text{and} \quad
t_i \ge \frac{i+2}{2},
\quad i=0,1,\cdots.
\label{eq:ti_sum}
\end{align}

Algorithm FGM1 is in Algorithm Class FO
\cite[Proposition 2]{drori:14:pof} with:
\begin{align}
\bhkip = \begin{cases}
	\frac{t_i - 1}{t_{i+1}}\bhki, & k = 0, \cdots, i - 2, \\
	\frac{t_i - 1}{t_{i+1}}(\bhimi - 1), & k = i - 1, \\
	1 + \frac{t_i - 1}{t_{i+1}}, & k = i,
\end{cases}
\label{eq:fgm1,h} 
\end{align}
for $i=0,\cdots,N-1$. 
Note that Algorithm FO with~\eqref{eq:fgm1,h} is impractical
as written for large-scale optimization problems,
whereas the mathematically equivalent version FGM1
is far more useful practically due to its efficient form.

While the sequence $\{\bmx_0,\cdots,\bmx_{N-1},\bmy_N\}$ of FGM1 
can be also written in class FO~\cite[Proposition 2]{drori:14:pof},
only the \emph{\auxiliary}~sequence $\{\bmy_0,\cdots,\bmy_N\}$ is 
known to achieve the rate $O\paren{\frac{1}{N^2}}$
for decreasing $f$
\cite{beck:09:afi,nesterov:83:amo}.
\DT~conjectured that
the \emph{\primary} sequence $\{\bmx_0,\cdots,\bmx_N\}$ of FGM1
also achieves the same $O\paren{\frac{1}{N^2}}$ rate
based on the numerical results 
using the PEP approach~\cite[Conjecture 2]{drori:14:pof};
our Section~\ref{sec:fgm,pep,ana} verifies the conjecture
by providing an analytical bound using the PEP approach.

\subsection{Nesterov's fast gradient method 2 (FGM2)}

In~\cite{nesterov:05:smo}, Nesterov proposed another fast gradient method
that has a different\footnote{
The fast gradient method in~\cite{nesterov:05:smo}
was originally developed to generalize FGM1 to the constrained case.
Here, this second form is introduced for use in later proofs.
}
form than FGM1
and that used a choice of $t_i$ factors different from~\eqref{eq:ti}.
Here, we use~\eqref{eq:ti} because it leads to faster convergence
than the factors used in~\cite{nesterov:05:smo}.
The algorithm in~\cite{nesterov:05:smo} then becomes FGM2 shown below.

\fbox{
\begin{minipage}[t]{0.85\linewidth}
\vspace{-10pt}
\begin{flalign*}
&\quad \text{\bf Algorithm FGM2} & \\
&\qquad \text{Input: } f\in C_L^{1,1}(\Reals^d)\text{ convex},\; \bmx_0\in\Reals^d,\; 
	t_0 = 1. & \\
&\qquad \text{For } i = 0,\cdots,N-1 & \\
&\qquad \qquad \bmy_{i+1} = \bmx_i - \frac{1}{L}f'(\bmx_i) & \\
&\qquad \qquad \bmz_{i+1} = \bmx_0 - \frac{1}{L}\sum_{k=0}^i t_k f'(\bmx_k) & \\
&\qquad \qquad t_{i+1} = \frac{1+\sqrt{1+4t_i^2}}{2} & \\
&\qquad \qquad \bmx_{i+1} = \paren{1 - \frac{1}{t_{i+1}}}\bmy_{i+1} 
				+ \frac{1}{t_{i+1}}\bmz_{i+1} & 
\end{flalign*}
\end{minipage}
} \vspace{5pt}

\noindent
Similar to FGM1,
the following proposition shows that FGM2 is in class FO with
\begin{align}
\bhkip = \begin{cases}
	\frac{1}{t_{i+1}}\paren{t_k - \sum_{j=k+1}^i \bhkj}, 
		& k = 0,\cdots,i - 1, \\
	1 + \frac{t_i - 1}{t_{i+1}}, & k = i,
\label{eq:fgm2,h} 
\end{cases}
\end{align}
for $i=0,\cdots,N-1$ with $t_i$ in~\eqref{eq:ti}.

\begin{proposition}
The sequence $\{\bmx_0,\cdots,\bmx_N\}$ generated by Algorithm FO 
with~\eqref{eq:fgm2,h} is identical
to the corresponding sequence generated by Algorithm FGM2.
\begin{proof}
We use induction,
and for clarity, we use the notation $\bmx_0',\cdots,\bmx_N'$
for Algorithm FO.
Clearly $\bmx_0' = \bmx_0$.
To prove equivalence for $i=1$:
\begin{align*}
\bmx_1' &= \bmx_0' - \frac{1}{L}\bhzo f'(\bmx_0') 
	= \bmx_0 - \frac{1}{L}\paren{1 + \frac{t_0 - 1}{t_1}}f'(\bmx_0) \\
	&= \paren{1 - \frac{1}{t_1} + \frac{1}{t_1}}
		\paren{\bmx_0 - \frac{1}{L}f'(\bmx_0)} 
	= \paren{1-\frac{1}{t_1}}\bmy_1 + \frac{1}{t_1}\bmz_1 = \bmx_1
.\end{align*}
Assuming $\bmx_i' = \bmx_i$ for $i = 0,\cdots,n$, we then have
{\allowdisplaybreaks[4]
\begin{align*}
\bmx_{n+1}' 
=& \bmx_n' - \frac{1}{L}\bhnnp f'(\bmx_n') 
	- \frac{1}{L}\sum_{k=0}^{n-1}\bhknp f'(\bmx_k') \\
=& \bmx_n - \frac{1}{L}\paren{1 + \frac{t_n - 1}{t_{n+1}}}f'(\bmx_n)
	- \frac{1}{L}\sum_{k=0}^{n-1}\frac{1}{t_{n+1}}
	\paren{t_k - \sum_{j=k+1}^n \bhkj} f'(\bmx_k) \\
=& \paren{1 - \frac{1}{t_{n+1}}}\paren{\bmx_n - \frac{1}{L}f'(\bmx_n)} \\
&	+ \frac{1}{t_{n+1}}\paren{\bmx_n 
	+ \frac{1}{L}\sum_{k=0}^{n-1}\sum_{j=k+1}^n \bhkj f'(\bmx_k) 
	- \frac{1}{L}\sum_{k=0}^n t_k f'(\bmx_k)} \\
=& \paren{1 - \frac{1}{t_{n+1}}} \bmy_{n+1}      
        + \frac{1}{t_{n+1}}\paren{\bmx_n         
	+ \frac{1}{L}\sum_{j=1}^n\sum_{k=0}^{j-1} \bhkj f'(\bmx_k)
        - \frac{1}{L}\sum_{k=0}^n t_k f'(\bmx_k)} \\
=& \paren{1 - \frac{1}{t_{n+1}}} \bmy_{n+1}                 
        + \frac{1}{t_{n+1}}\paren{\bmx_0
        - \frac{1}{L}\sum_{k=0}^n t_k f'(\bmx_k)} 
= \bmx_{n+1}
.\end{align*}}
$\!\!\!$The fifth equality uses the telescoping sum
$\bmx_n = \bmx_0 + \sum_{j=1}^{n} (\bmx_j - \bmx_{j-1})$
and~\eqref{eq:fo} in Algorithm FO.
\qed
\end{proof}
\end{proposition}

We show next the equivalence of Nesterov's two algorithms FGM1 and FGM2
for smooth unconstrained convex minimization
using~\eqref{eq:fgm1,h} and~\eqref{eq:fgm2,h}.

\begin{proposition}
\label{prop:fgm1,fgm2}
The sequence $\{\bmx_0,\cdots,\bmx_N\}$ generated by Algorithm FGM2 is identical
to the corresponding sequence generated by Algorithm FGM1.
\begin{proof}
We prove the statement by showing the equivalence of
\eqref{eq:fgm1,h} and~\eqref{eq:fgm2,h}.
We use the notation $\bhki'$ for the coefficients~\eqref{eq:fgm2,h} 
of Algorithm FGM2
to distinguish from those of Algorithm FGM1.

It is obvious that 
$\bhiip' = \bhiip,\; i=0,\cdots,N-1$,
and we can easily prove for $i=0,\cdots,N-1$ that
\begin{align*}
\bhimip' &= \frac{1}{t_{i+1}}\paren{t_{i-1} - \bhimi'}
        = \frac{1}{t_{i+1}}\paren{t_{i-1} - \paren{1 + \frac{t_{i-1} - 1}{t_i}}} \\
        &= \frac{(t_i - 1)(t_{i-1} - 1)}{t_i t_{i+1}}
        = \frac{t_i - 1}{t_{i+1}}\paren{\bhimi - 1} = \bhimip
.\end{align*}

We next use induction
by assuming $\bhkip' = \bhkip$ 
for $i=0,\cdots,n-1,\;k=0,\cdots,i$. We then have
\begin{align*}
\bhknp' &= \frac{1}{t_{n+1}}\paren{t_k - \sum_{j=k+1}^n \bhkj'}
	= \frac{1}{t_{n+1}}\paren{t_k - \sum_{j=k+1}^{n-1} \bhkj' - \bhkn'} \\
	&= \frac{t_n - 1}{t_{n+1}}\bhkn' = \frac{t_n - 1}{t_{n+1}}\bhkn = \bhknp 
\end{align*}
for $k=0,\cdots,n-2$.
Note that this proof is independent of the choice of $t_i$.
\qed
\end{proof}
\end{proposition}

\subsection{A convergence bound for Nesterov's fast gradient methods} 
\label{sec:conv,fgm}

Algorithms FGM1 and FGM2
generate the same sequences $\{\bmx_i\}$ and $\{\bmy_i\}$, 
and the \auxiliary~sequence $\{\bmy_i\}$ is known to satisfy the bound\footnote{
The second inequality of~\eqref{eq:fgm,conv} is widely known
since it provides simpler interpretation of a convergence bound,
compared to the first inequality of~\eqref{eq:fgm,conv}.
}
\cite{beck:09:afi,nesterov:83:amo,nesterov:05:smo}:
\begin{align}
f(\bmy_n) - f(\bmx_*) \le \frac{L||\bmx_0 - \bmx_*||^2}{2t_{n-1}^2}
\le \frac{2L||\bmx_0 - \bmx_*||^2}{(n+1)^2},
\quad \forall \bmx_* \in X_*(f)
\label{eq:fgm,conv}
\end{align}
for $n\geq1$,
which was the previously best known analytical bound
of first-order methods for smooth unconstrained convex minimization;
\DT's PEP approach provides
a tighter \emph{numerical} bound for the sequences $\{\bmx_i\}$ and $\{\bmy_i\}$ 
compared to the analytical bound~\eqref{eq:fgm,conv}
\cite[Table 1]{drori:14:pof}.
Using the PEP approach, Section~\ref{sec:fgm,pep,ana} provides 
a new \emph{analytical} bound
for the \primary~sequence $\{\bmx_i\}$ of FGM1 and FGM2.

Nesterov described a convex function $f \in C_L^{1,1}(\Reals^d)$
for which any first-order algorithm generating the sequence $\{\bmx_i\}$
in the class of Algorithm FO satisfies
\cite[Theorem 2.1.7]{nesterov:04}:
\begin{align}
\frac{3L||\bmx_0 - \bmx_*||^2}{32(n+1)^2} \le f(\bmx_n) - f(\bmx_*),
\quad \forall \bmx_* \in X_*(f)
\label{eq:nes,conv}
\end{align}
for $n = 1,\cdots,\floor{\frac{d-1}{2}}$,
indicating that Nesterov's two FGM1 and FMG2
achieve the optimal rate $O\paren{\frac{1}{N^2}}$.
(Note that the bound~\eqref{eq:nes,conv}
is valid if the large-scale condition ``$d\ge 2N+1$'' is satisfied.)
However,~\eqref{eq:nes,conv} also illustrates the potential room for improving
first-order algorithms by a constant factor.

To narrow this gap, \DT~\cite{drori:14:pof}
used a relaxation of problem~\eqref{eq:HP} to find 
the ``optimal'' choice of $\{\hki\}$ for Algorithm FO
that minimizes a relaxed bound on $f(\bmx_N) - f(\bmx_*)$ at the $N$th iteration,
which was found numerically to provide a twice smaller bound than~\eqref{eq:fgm,conv},
yet remained computationally impractical.

We next review the PEP approach for solving a relaxed version of~\eqref{eq:PEP}.

\section{\DT's convergence bound for first-order algorithms using PEP}
\label{sec:pep}

This section summarizes the relaxation scheme for the PEP approach
that transforms problem~\eqref{eq:PEP} into a tractable form~\cite{drori:14:pof}.
The relaxed PEP bounds are used in later sections.

Problem~\eqref{eq:PEP} is challenging to solve
due to the (infinite-dimensional) functional constraint on $f$,
so \DT~\cite{drori:14:pof} cleverly relax the constraint by using
a well-known property for the class of convex $C_L^{1,1}$ functions
in~\cite[Theorem 2.1.5]{nesterov:04} 
and further relax as follows:
\begin{align*}
\mathcal{B}_{\mathrm{P1}}(\bmh,N,\dd,L,R) \defequ\;
&\max_{\substack{\bmG\in\Reals^{(N+1)d}, \\ \bmdel\in\Reals^{N+1}}}
        LR^2\delta_N 
	\label{eq:pPEP} \tag{P1} \\
	&\st \; \Tr{\bmG^\top\bmA_{i-1,i}(\bmh)\bmG} \le \delta_{i-1} - \delta_i,
		\quad i=1,\cdots,N, \nonumber \\
	&\quad\;\;\; \Tr{\bmG^\top\bmD_i(\bmh)\bmG + \bmnu \bmu_i^\top\bmG} \le - \delta_i,   
                \quad i = 0,\cdots,N, \nonumber 
\end{align*}
for any given unit vector $\bmnu\in\Reals^d$,
by defining
$\delta_i \defequ \frac{1}{L||\bmx_0 - \bmx_*||^2}(f(\bmx_i) - f(\bmx_*))$ and 
$\bmg_i \defequ \frac{1}{L||\bmx_0 - \bmx_*||}f'(\bmx_i)$
for $i=0,\cdots,N,*$,
and denoting
the unit vectors\footnote{
The vector $\bme_{N,i}^{ }$ 
is the $i$th standard basis vector in $\Reals^{N}$,
having $1$ for the $i$th entry
and zero for all other elements.}
$\bmu_i = \bme_{N+1,i+1}^{ } \in \Reals^{N+1}$, 
the $(N+1)\times 1$ vector $\bmdel = [\delta_0,\;\cdots,\;\delta_N]^\top$,
the $(N+1)\times d$ matrix $\bmG = [\bmg_0,\;\cdots,\;\bmg_N]^\top$, 
and the $(N+1)\times(N+1)$ symmetric matrices:
{\allowdisplaybreaks[4]
\begin{align}
\begin{cases}
\bmA_{i-1,i}(\bmh) \defequ \frac{1}{2}(\bmu_{i-1} - \bmu_i)(\bmu_{i-1} - \bmu_i)^\top
        + \frac{1}{2}\sum_{k=0}^{i-1} 
        \hki (\bmu_i\bmu_k^\top + \bmu_k\bmu_i^\top), & \\
\bmD_i(\bmh) \defequ \frac{1}{2}\bmu_i\bmu_i^\top + \frac{1}{2}\sum_{j=1}^i\sum_{k=0}^{j-1}
                \hkj(\bmu_i\bmu_k^\top + \bmu_k\bmu_i^\top). &
\end{cases}
\label{eq:ABCD}
\end{align}
}

\DT~\cite{drori:14:pof} finally use a duality approach
on~\eqref{eq:pPEP}.
Replacing $\max_{\bmG,\bmdel} LR^2\delta_N$
by $\min_{\bmG,\bmdel} -\delta_N$ for convenience,
the Lagrangian of the corresponding 
constrained minimization problem~\eqref{eq:pPEP} becomes
the following separable function in $(\bmdel,\bmG)$:
\begin{align*}
\cL(\bmG,\bmdel,\bmlam,\bmtau;\bmh) 
	= \cL_1(\bmdel,\bmlam,\bmtau) + \cL_2(\bmG,\bmlam,\bmtau;\bmh) \nonumber
,\end{align*}
where  
\begin{align*}	
\begin{cases}
\cL_1(\bmdel,\bmlam,\bmtau) 
	\defequ - \delta_N + \sum_{i=1}^N \lambda_i(\delta_i - \delta_{i-1})
		+ \sum_{i=0}^N \tau_i\delta_i, & \\
\cL_2(\bmG,\bmlam,\bmtau;\bmh)
	\defequ \sum_{i=1}^N\lambda_i\Tr{\bmG^\top\bmA_{i-1,i}(\bmh)\bmG}
		+ \sum_{i=0}^N\tau_i\Tr{\bmG^\top\bmD_i(\bmh)\bmG + \bmnu\bmu_i^\top\bmG}, & 
\end{cases}
\end{align*}
with dual variables $\bmlam = (\lambda_1,\cdots,\lambda_N)^\top\in\Reals_+^N$
and $\bmtau = (\tau_0,\cdots,\tau_N)^\top\in\Reals_+^{N+1}$.
The corresponding dual function is defined as 
\begin{align}
H(\bmlam,\bmtau;\bmh) 
		&= \min_{\bmdel\in\Reals^{N+1}} \cL_1(\bmdel,\bmlam,\bmtau) 
		+ \min_{\bmG\in\Reals^{(N+1)d}} \cL_2(\bmG,\bmlam,\bmtau;\bmh)
\label{eq:H}
.\end{align}
Here $\min_{\bmdel} \cL_1(\bmdel,\bmlam,\bmtau) = 0$ 
for any $(\bmlam,\bmtau)\in\Lambda$,
where 
\begin{align}
\Lambda = \bigg\{
	(\bmlam,\bmtau)\in\Reals_+^N\times\Reals_+^{N+1}
        \;:\; \begin{array}{l}
		\tau_0 = \lambda_1,\; \lambda_N + \tau_N = 1 \\
                \lambda_i - \lambda_{i+1} + \tau_i = 0,\;  
                 i=1,\cdots,N-1 
		\end{array}\bigg\},        
	\label{eq:Lambda}
\end{align}
and $-\infty$ otherwise. 
In~\cite{drori:14:pof},
the dual function~\eqref{eq:H} for any given unit vector $\bmnu\in\Reals^d$  
was found to be 
\begin{align}
H(\bmlam,\bmtau;\bmh) &= \min_{\bmw\in\Reals^{N+1}} 
		\braces{\bmw^\top\bmS(\bmh,\bmlam,\bmtau)\bmw + \bmtau^\top\bmw} 
		\nonumber \\ 
	&= \max_{\gamma\in\Reals}
	        \braces{- \frac{1}{2}\gamma
		\;:\; \bmw^\top\bmS(\bmh,\bmlam,\bmtau)\bmw + \bmtau^\top\bmw 
		\ge - \frac{1}{2}\gamma,\; \forall \bmw\in\Reals^{N+1}} 
		\nonumber \\
	&= \max_{\gamma\in\Reals} \braces{- \frac{1}{2}\gamma
		\;:\; \paren{\begin{array}{cc}
        	\bmS(\bmh,\bmlam,\bmtau) & \frac{1}{2}\bmtau \\
        	\frac{1}{2}\bmtau^\top & \frac{1}{2}\gamma
        	\end{array}} \succeq 0}
		\label{eq:prevH}
\end{align}
for any given $(\bmlam,\bmtau) \in \Lambda$,
where \DT~\cite{drori:14:pof} define
the following $(N+1)\times(N+1)$ matrix using
the definition of $\bmA_{i-1,i}(\bmh)$ and $\bmD_i(\bmh)$ in~\eqref{eq:ABCD}:
\begin{align}
\bmS(\bmh,\bmlam,\bmtau) \defequ& \sum_{i=1}^N \lambda_i\bmA_{i-1,i}(\bmh) 
		+ \sum_{i=0}^N\tau_i\bmD_i(\bmh) \nonumber \\
=& \frac{1}{2}\sum_{i=1}^N \lambda_i(\bmu_{i-1} - \bmu_i)(\bmu_{i-1} - \bmu_i)^\top
	+ \frac{1}{2}\sum_{i=0}^N \tau_i\bmu_i\bmu_i^\top \nonumber\\
	&+ \frac{1}{2}\sum_{i=1}^N\sum_{k=0}^{i-1} 
	\paren{\lambda_i\hki + \tau_i\sum_{j=k+1}^i\hkj}
	(\bmu_i\bmu_k^\top + \bmu_k\bmu_i^\top).
\label{eq:S}
\end{align}
In short, using the dual approach 
on the problem~\eqref{eq:pPEP} 
yields the following bound:
\begin{align}
\mathcal{B}_{\mathrm{D}}(\bmh,N,L,R) \defequ
\min_{\substack{\bmlam\in\Reals^N, \\ \bmtau\in\Reals^{N+1}, \\ \gamma\in\Reals}} 
	\braces{ \frac{1}{2}LR^2\gamma\;:\;
	\paren{\begin{array}{cc}
	\bmS(\bmh,\bmlam,\bmtau) & \frac{1}{2}\bmtau \\
	\frac{1}{2}\bmtau^\top & \frac{1}{2}\gamma
	\end{array}} \succeq 0,
	\quad
	(\bmlam,\bmtau) \in \Lambda},
	\label{eq:dPEP} \tag{D}
\end{align}
recalling that we previously replaced $\max_{\bmG,\bmlam} LR^2\delta_N$
by $\min_{\bmG,\bmlam} -\delta_N$ for convenience.
Problem~\eqref{eq:dPEP} can be solved
using any numerical SDP method~\cite{cvx,gb08}
for given \bmh and $N$,
noting that $R$ is just a multiplicative scalar in~\eqref{eq:dPEP}.
Interestingly, this bound $\mathcal{B}_{\mathrm{D}}(\bmh,N,L,R)$
is independent of dimension $d$.

Overall, 
\DT~\cite{drori:14:pof} introduced a series of relaxations to
the problem~\eqref{eq:PEP},
eventually reaching the solvable problem~\eqref{eq:dPEP} that
provides a valid upper bound as
\begin{align*}
f(\bmx_N) - f(\bmx_*) \le
\mathcal{B}_{\mathrm{P}}(\bmh,N,\dd,L,R)
\le \mathcal{B}_{\mathrm{D}}(\bmh,N,L,R)
\end{align*}
where
$\bmx_N$ is generated by Algorithm FO with given \bmh and $N$,
and $||\bmx_0 - \bmx_*|| \le R$.
This bound is for a given \bmh
and later we optimize the bound over \bmh.

Solving problem~\eqref{eq:dPEP} with a SDP method 
for any given coefficients \bmh and $N$
provides
a numerical convergence bound for $f(\bmx_N) - f(\bmx_*)$~\cite{drori:14:pof}.
However, numerical bounds only partially explain the behavior of
algorithms in class FO.
An analytical bound 
of gradient methods with a constant step $0<h\le1$,
for example,
was found using a specific PEP approach~\cite{drori:14:pof},
but no other analytical bound was discussed in~\cite{drori:14:pof}.
The next section exploits the PEP approach
to reveal a new analytical bound 
for the \primary~sequence $\{f(\bmx_i)\}$ generated by FGM1 or FGM2 as an example,
confirming the conjecture by \DT~that
the \emph{\primary}~sequence $\{\bmx_i\}$ achieves the same rate $O\paren{\frac{1}{N^2}}$
as the \emph{\auxiliary}~sequence $\{\bmy_i\}$~\cite[Conjecture 2]{drori:14:pof}.

\section{A new analytical bound for Nesterov's fast gradient methods}
\label{sec:fgm,pep,ana}

This section provides an analytical bound
for the \primary~sequence $\{\bmx_i\}$ in FGM1 and FGM2.

For the \bh factors in~\eqref{eq:fgm1,h} or~\eqref{eq:fgm2,h}
of Nesterov's fast gradient methods,
the following choice of dual variables
(inspired by Section~\ref{sec:opt,fo,ana})
is a feasible point of problem~\eqref{eq:dPEP}:
\begin{align}
\blam_i &= \frac{t_{i-1}^2}{t_N^2},\quad i=1,\cdots,N,
	\label{eq:blam} \\
\btau_i &= \frac{t_i}{t_N^2},\quad i=0,\cdots,N, 
	\label{eq:btau} \\
\bgam &= \frac{1}{t_N^2}, \label{eq:bgam}
\end{align}
with $t_i$ in~\eqref{eq:ti},
as shown in the following lemma.

\begin{lemma}
\label{lem:fgm,fea}
The choice $(\bmblam,\bmbtau,\bgam)$ 
in~\eqref{eq:blam},~\eqref{eq:btau} and~\eqref{eq:bgam}
is a feasible point of the problem~\eqref{eq:dPEP} 
for the \bh designs given in~\eqref{eq:fgm1,h} or~\eqref{eq:fgm2,h}
that are used in Nesterov's FGM1 and FGM2.
\begin{proof}
It is obvious that $(\bmblam,\bmbtau)\in\Lambda$ 
using $t_i^2 = \sum_{k=0}^it_k$ in~\eqref{eq:ti_sum}.
We next rewrite $\bmS(\bh,\bmblam,\bmbtau)$
using~\eqref{eq:fgm2,h},~\eqref{eq:blam} and~\eqref{eq:btau}
to show that
the choice $(\bmblam,\bmbtau,\bgam)$
satisfies the positive semidefinite condition in~\eqref{eq:dPEP}
for given \bh.

For any \bmh and $(\bmlam,\bmtau) \in \Lambda$,
the $(i,k)$-th entry of 
the symmetric matrix $\bmS(\bmh,\bmlam,\bmtau)$ in~\eqref{eq:S}
can be written as
\begin{align}
S_{i,k}(\bmh,\bmlam,\bmtau)
&= \begin{cases}
	\frac{1}{2}\Big((\lambda_{i} + \tau_{i})\hki
                + \tau_{i}\sum_{j=k+1}^{i-1} &\!\!\!\!\!\!\hkj\Big), \\
        & i=2,\cdots,N,\;k =0,\cdots, i - 2, \\
        \frac{1}{2}\paren{(\lambda_{i} + \tau_{i})\hki - \lambda_i},
        & i=1,\cdots,N,\;k = i-1, \\
        \lambda_{i+1},
        & i=0,\cdots,N-1,\;k = i, \\
        \frac{1}{2},
        & i=N,\;k = i. 
	\end{cases}
\label{eq:S2}
\end{align}
Inserting \bh~\eqref{eq:fgm2,h}, \bmblam~\eqref{eq:blam} 
and \bmbtau~\eqref{eq:btau}
into~\eqref{eq:S2}
and using $\blam_i + \btau_i = \frac{t_i^2}{t_N^2}$ for $i=1,\cdots,N$, we get
{\allowdisplaybreaks[4]
\begin{align*}
S_{i,k}(\bh,\bmblam,\bmbtau)
&= \begin{cases}
        \frac{1}{2}\Big(\frac{t_{i}^2}{t_N^2}
                        \frac{1}{t_{i}}\paren{t_k - \sum_{j=k+1}^{i-1}\bhkj} 
                        +&\!\!\!\!\!\! \frac{t_{i}}{t_N^2}\sum_{j=k+1}^{i-1}\bhkj\Big), \\
	& i=2,\cdots,N,\;k =0,\cdots,i - 2, \\
        \frac{1}{2}\paren{\frac{t_{i}^2}{t_N^2}
                \paren{1 + \frac{t_{i-1} - 1}{t_{i}}} - \frac{t_{i-1}^2}{t_N^2}},
	& i=1,\cdots,N,\;k = i-1, \\
        \frac{t_i^2}{t_N^2},
        & i=0,\cdots,N-1,\;k = i, \\
        \frac{1}{2},
        & i=N,\;k = i, 
        \end{cases} \\
&= \begin{cases}
	\frac{t_it_k}{2t_N^2}
        & i=1,\cdots,N,\;k =0,\cdots, i - 1, \\
        \frac{t_i^2}{t_N^2},
        & i=0,\cdots,N-1,\;k = i, \\
        \frac{t_N^2}{2t_N^2},
        & i=N,\;k = i,
        \end{cases} \\
&= \frac{1}{2t_N^2}\paren{\bmt\,\bmt^\top + \diag{(\bmtt^\top,0)}},
\end{align*}}
$\!\!\!$where $\bmt = (t_0, \cdots, t_N)^\top$
and $\bmtt = (t_0^2, \cdots, t_{N-1}^2)^\top$.
The second equality uses $t_i^2 - t_i - t_{i-1}^2 = 0$ in~\eqref{eq:ti_sum},
and $\diag{\bmt}$ denotes a matrix where diagonal elements
are filled with elements of a vector \bmt
and zero for other elements.

Finally, using \bgam in~\eqref{eq:bgam}, we have
\begin{align*}
\paren{\begin{array}{cc}
        \bmS(\bh,\bmblam,\bmbtau) & \frac{1}{2}\bmbtau \\
        \frac{1}{2}\bmbtau^\top & \frac{1}{2}\bgam
\end{array}}
&= \paren{\begin{array}{cc} \frac{1}{2t_N^2}\paren{\bmt\,\bmt^\top + \diag{(\bmtt^\top,0)}} 
		& \frac{1}{2t_N^2}\bmt \\ 
		\frac{1}{2t_N^2}\bmt^\top 
		& \frac{1}{2t_N^2} \end{array}} \\
&= \frac{1}{2t_N^2}\braces{\paren{\begin{array}{c} \bmt \\ 1 \end{array}}
		\paren{\begin{array}{c} \bmt \\ 1 \end{array}}^\top 
		+ \diag{(\bmtt^\top,0,0)}}\succeq 0
.\end{align*}
\qed
\end{proof}
\end{lemma}

Using Lemma~\ref{lem:fgm,fea},
we provide an analytical convergence bound for
the \primary~sequence $\{\bmx_i\}$ of FGM1 and FMG2.

\begin{theorem}
\label{thm:fgm,aux}
Let $f\;:\;\Reals^d\rightarrow\Reals$ be convex and $C_L^{1,1}$,
and let $\bmx_0,\bmx_1,\cdots \in \Reals^d$
be generated by FGM1 or FGM2.
Then for $n\ge1$,
\begin{align}
f(\bmx_n) - f(\bmx_*) \le \frac{L||\bmx_0 - \bmx_*||^2}{2t_n^2}
\le \frac{2L||\bmx_0 - \bmx_*||^2}{(n+2)^2},
\quad \forall \bmx_* \in X_*(f)
\label{eq:fgm,conv,aux}
.\end{align}
\begin{proof}
Using \bgam~\eqref{eq:bgam} 
and $t_N^2 \ge \frac{(N+2)^2}{4}$ from~\eqref{eq:ti_sum}, we have
\begin{align}
f(\bmx_N) - f(\bmx_*) \le 
		\mathcal{B}_{\mathrm{D}}(\bh,N,L,R)
		\le \frac{1}{2}LR^2\bgam \le \frac{2LR^2}{(N+2)^2},
\quad \forall \bmx_* \in X_*(f)
\label{eq:fgm,conv,aux1}
\end{align}
for given \bh in~\eqref{eq:fgm1,h} or~\eqref{eq:fgm2,h},
based on Lemma~\ref{lem:fgm,fea}.
Since the coefficients \bh in~\eqref{eq:fgm1,h} or~\eqref{eq:fgm2,h}
are recursive and
do not depend on a given $N$,
we can extend~\eqref{eq:fgm,conv,aux1} for all iterations ($n\ge1$).
Finally, we let $R = ||\bmx_0 - \bmx_*||$.
\qed
\end{proof}
\end{theorem}

Theorem~\ref{thm:fgm,aux} illustrates using the PEP approach
to find an analytical bound for an algorithm in class FO.
We used a SDP solver~\cite{cvx,gb08}
to verify numerically that the choice $(\bmblam,\bmbtau,\bgam)$
in~\eqref{eq:blam},~\eqref{eq:btau} and~\eqref{eq:bgam}
is not an optimal solution of~\eqref{eq:dPEP}
for given \bh in~\eqref{eq:fgm1,h} or~\eqref{eq:fgm2,h}.
Nevertheless, this feasible point $(\bmblam,\bmbtau,\bgam)$
provides a valid upper bound for the sequence $\{\bmx_i\}$
of FGM1 and FGM2
as shown in Theorem~\ref{thm:fgm,aux}
that is similar to~\eqref{eq:fgm,conv}
and verifies \DT's conjecture~\cite[Conjecture 2]{drori:14:pof}.

The next section reviews \DT's work~\cite{drori:14:pof}
on numerically optimizing step coefficients \bmh
in the class of first-order methods
over the relaxed convergence bound $\mathcal{B}_{\mathrm{D}}(\bmh,N,L,R)$.
Then, we find an analytical form of the optimized step coefficients
and explicitly show that Algorithm FO with such coefficients 
achieves a convergence bound
that is twice as small as
\eqref{eq:fgm,conv} and~\eqref{eq:fgm,conv,aux}.

\section{Towards optimized first-order algorithms}
\label{sec:opt,fo}

\subsection{\DT's numerically optimized first-order algorithms}
\label{sec:opt,fo,num}

This section summarizes the numerically optimized first-order algorithms 
described in~\cite{drori:14:pof}.

Having relaxed~\eqref{eq:PEP} in Section~\ref{sec:pep}
to~\eqref{eq:dPEP},
\DT~proposed to optimize \bmh 
by relaxing~\eqref{eq:HP} as follows:
\begin{align}
\Hh \defequ
\argmin{\bmh\in\Reals^{\Frac{N(N+1)}{2}}}\;
\mathcal{B}_{\mathrm{D}}(\bmh,N,L,R)
        \label{eq:HD} \tag{HD}
,\end{align}
where \Hh is independent of both $L$ and $R$,
since $\mathcal{B}_{\mathrm{D}}(\bmh,N,L,R) = LR^2\mathcal{B}_{\mathrm{D}}(\bmh,N,1,1)$.
Problem~\eqref{eq:HD} is a bilinear optimization problem
in terms of \bmh and the dual variables in~\eqref{eq:dPEP},
unlike the linear SDP problem~\eqref{eq:dPEP}.
To simplify,
\DT~\cite{drori:14:pof}
introduced a variable $\bmr = \{\rki\}_{0\le k < i \le N}$:
\begin{align}
\rki = \lambda_i \hki + \tau_i\sum_{j=k+1}^i\hkj
\label{eq:rki}
\end{align}
to convert~\eqref{eq:HD} into
the following linear SDP problem:
\begin{align}
\Hr \defequ
\argmin{\bmr\in\Reals^{\Frac{N(N+1)}{2}}}\;
\breve{\mathcal{B}}_{\mathrm{D}}(\bmr,N,L,R),
\label{eq:HDR} \tag{RD}
\end{align}
where
\begin{align}
\breve{\mathcal{B}}_{\mathrm{D}}(\bmr,N,L,R) \defequ&
\min_{\substack{\bmlam\in\Reals^N, \\ \bmtau\in\Reals^{N+1}, \\ \gamma\in\Reals}}
        \braces{ \frac{1}{2}LR^2\gamma\;:\;
        \paren{\begin{array}{cc}
        \Sr(\bmr,\bmlam,\bmtau) & \frac{1}{2}\bmtau \\
        \frac{1}{2}\bmtau^\top & \frac{1}{2}\gamma
        \end{array}} \succeq 0,
        (\bmlam,\bmtau) \in \Lambda}, \nonumber \\
\Sr(\bmr,\bmlam,\bmtau) 
\defequ& \frac{1}{2}\sum_{i=1}^N \lambda_i(\bmu_{i-1} - \bmu_i)(\bmu_{i-1} - \bmu_i)^\top
        + \frac{1}{2}\sum_{i=0}^N \tau_i\bmu_i\bmu_i^\top \nonumber\\
        &+ \frac{1}{2}\sum_{i=1}^N\sum_{k=0}^{i-1}
        \rki(\bmu_i\bmu_k^\top + \bmu_k\bmu_i^\top).
\label{eq:Sr}
\end{align}
An optimal solution $(\Hr,\bmhlam,\bmhtau,\hgam)$ of~\eqref{eq:HDR}
for a given $N$
can be computed by any numerical SDP method~\cite{cvx,gb08}.
\DT~showed that
the resulting values $(\bmhlam,\bmhtau,\hgam)$
with the following \Hh:
\begin{align}
\Hhki = \begin{cases}
\frac{\Hrki - \htau_i\sum_{j=k+1}^{i-1}\Hhkj}{\hlam_i + \htau_i},
& \hlam_i + \htau_i \neq 0, \\
0, & \text{otherwise},
\end{cases}
\label{eq:hki,rki}
\end{align}
for $i=1,\cdots,N,\;k=0,\cdots,i-1$
become an optimal solution of~\eqref{eq:HD}
\cite[Theorem 3]{drori:14:pof},\footnote{
Equation (5.2) in~\cite[Theorem 3]{drori:14:pof} 
that is derived from~\eqref{eq:rki} has typos
that we fixed in~\eqref{eq:hki,rki}.
}
where both~\eqref{eq:HD} and~\eqref{eq:HDR}
achieve the same optimal value,
\ie, $\mathcal{B}_{\mathrm{D}}(\Hh,N,L,R) = \breve{\mathcal{B}}_{\mathrm{D}}(\Hr,N,L,R)$.

The numerical results for problem~\eqref{eq:HD} in~\cite{drori:14:pof} 
provided a convergence bound that is about two-times smaller  
than that of Nesterov's fast gradient methods
for a couple of choices of $N$ in~\cite[Tables 1 and 2]{drori:14:pof}.
However, numerical calculations cannot verify the acceleration for all $N$,
and SDP computation for solving~\eqref{eq:HDR} becomes expensive for large $N$.
In the next section, we analytically solve problem~\eqref{eq:HD},
which is our first main contribution.

\subsection{Proposed analytically optimized first-order algorithms}
\label{sec:opt,fo,ana}

This section provides an analytical optimal solution of~\eqref{eq:HD} 
by reformulating~\eqref{eq:HDR}.

We first find an equivalent form of 
the dual function $H(\bmlam,\bmtau;\bmh)$ in~\eqref{eq:H}
that differs from~\eqref{eq:prevH}
by using the following equality:
\begin{align} 
S_{N,N}(\bmh,\bmlam,\bmtau) = \frac{1}{2}
\text{ for any } (\bmlam,\bmtau) \in \Lambda
\label{eq:fact}
,\end{align}
\ie, the $(N,N)$-th entry of $\bmS(\bmh,\bmlam,\bmtau)$ 
in~\eqref{eq:S} and~\eqref{eq:S2}
is $\frac{1}{2}$ for any $(\bmlam,\bmtau) \in \Lambda$.
Hereafter we use the notation
\begin{align}
\bmS(\bmh,\bmlam,\bmtau) \defequ \paren{\begin{array}{cc}
		\tS(\bmh,\bmlam,\bmtau) & \ts(\bmh,\bmlam,\bmtau) \\
		\ts(\bmh,\bmlam,\bmtau)^\top & \frac{1}{2}
	\end{array}},
\quad \bmw \defequ \paren{\begin{array}{c} \tw \\ w_N \end{array}},
\quad \text{and} 
\quad \bmtau \defequ \paren{\begin{array}{c} \ttau \\ \tau_N \end{array}} 
\label{eq:notation}
,\end{align}
where $\tS(\bmh,\bmlam,\bmtau)$ is a $N\times N$ symmetric matrix,
$\ts(\bmh,\bmlam,\bmtau)$, $\tw$ and $\ttau$ are $N\times1$ vectors,
and $w_N$ and $\tau_N$ are scalars.
We omit the arguments $(\bmh,\bmlam,\bmtau)$ 
in $\tS(\bmh,\bmlam,\bmtau)$ and $\ts(\bmh,\bmlam,\bmtau)$
for notational simplicity in the next derivation.
For any given $(\bmlam,\bmtau) \in \Lambda$, 
we rewrite $H(\bmlam,\bmtau;\bmh)$ in~\eqref{eq:H} and~\eqref{eq:prevH} as follows:
\begin{align}
H(\bmlam,\bmtau;\bmh)  
=& \min_{\bmw\in\Reals^{N+1}}
	\braces{\tw^\top\tS\tw + \ttau^\top\tw
	+ 2\tw^\top\ts w_N + \frac{1}{2}w_N^2 + \tau_Nw_N}
	\nonumber \\
=& \min_{\tw\in\Reals^N}
	\braces{\tw^\top(\tS - 2\ts\ts^\top)\tw
	+ (\ttau - 2\ts\tau_N)^\top\tw - \frac{1}{2}\tau_N^2}
	\nonumber \\
=& \max_{\gamma\in\Reals}
        \braces{- \frac{1}{2}\gamma\;:\;
	\begin{array}{l}
	\tw^\top(\tS - 2\ts\ts^\top)\tw
        + (\ttau - 2\ts\tau_N)^\top\tw - \frac{1}{2}\tau_N^2
	\ge - \frac{1}{2}\gamma, \\
	\forall\tw\in\Reals^N 
        \end{array}}
	\nonumber \\
=& \max_{\gamma\in\Reals}
        \braces{- \frac{1}{2}\gamma\;:\;
	\paren{\begin{array}{cc}
        \tS - 2\ts\ts^\top & \frac{1}{2}(\ttau - 2\ts\tau_N) \\
	\frac{1}{2}(\ttau - 2\ts\tau_N)^\top & \frac{1}{2}(\gamma - \tau_N^2)
	\end{array}} \succeq 0} 
	\label{eq:newH}
,\end{align}
where the second  
equality comes from minimizing the function 
with respect to $w_N$.

Using~\eqref{eq:newH} instead of~\eqref{eq:prevH}
for the function $H(\bmlam,\bmtau;\bmh)$
and again using the variable \bmr in~\eqref{eq:rki}
leads to the following optimization problem
that is equivalent to~\eqref{eq:HDR}:
\begin{align}
\Hr = 
\argmin{\bmr\in\Reals^{\Frac{N(N+1)}{2}}}\;
\breve{\mathcal{B}}_{\mathrm{D1}}(\bmr,N,L,R)
\label{eq:HDR1} \tag{RD1}
,\end{align}
where
\begin{align}
\breve{\mathcal{B}}_{\mathrm{D1}}(\bmr,N,L,R) \defequ&
\min_{\substack{\bmlam\in\Reals^N, \\ \bmtau\in\Reals^{N+1}, \\ \gamma\in\Reals}}
        \braces{ \frac{1}{2}LR^2\gamma\;:\;
        \paren{\begin{array}{cc}
        \Qr - 2\qr\qr^\top & \frac{1}{2}(\ttau - 2\qr\tau_N) \\
        \frac{1}{2}(\ttau - 2\qr\tau_N)^\top & \frac{1}{2}(\gamma - \tau_N^2)
        \end{array}} \succeq 0,\;
        (\bmlam,\bmtau) \in \Lambda}, \nonumber \\
\Qr(\bmr,\bmlam,\bmtau) 
	=& \frac{1}{2}\sum_{i=1}^{N-1}
		\lambda_i(\tbmu_{i-1} - \tbmu_i)(\tbmu_{i-1} - \tbmu_i)^\top
	+ \frac{1}{2}\lambda_N\tbmu_{N-1}\tbmu_{N-1}^\top 
	\nonumber \\
	&+ \frac{1}{2}\sum_{i=0}^{N-1}\tau_i\tbmu_i\tbmu_i^\top
	+ \frac{1}{2}\sum_{i=1}^{N-1}\sum_{k=0}^{i-1}
		\rki(\tbmu_i\tbmu_k^\top + \tbmu_k\tbmu_i^\top), 
\label{eq:Qr} \\
\qr(\bmr,\bmlam,\bmtau) 
	=& \frac{1}{2}\sum_{k=0}^{N-2}\rkN\tbmu_k,  
	+ \frac{1}{2}(\rNmN - \lambda_N)\tbmu_{N-1}
\label{eq:qr}
\end{align}
for $\tbmu_i = \bme_{N,i+1}^{ } \in \Reals^N$.
We omit the arguments $(\bmr,\bmlam,\bmtau)$
in $\Qr(\bmr,\bmlam,\bmtau)$ and $\qr(\bmr,\bmlam,\bmtau)$
for notational simplicity.
Unlike~\eqref{eq:HDR},
we observe that the new equivalent form~\eqref{eq:HDR1}
has a point 
at the boundary of the positive semidefinite condition,
and we will show that
the point is indeed an optimal solution of
both~\eqref{eq:HDR} and~\eqref{eq:HDR1}.

\begin{lemma}
\label{lem:feas}
A feasible point of both~\eqref{eq:HDR} and~\eqref{eq:HDR1}
is $(\Hr,\bmhlam,\bmhtau,\hgam)$, where
\begin{align}
\Hrki &= \begin{cases}
\frac{4\theta_i\theta_k}{\theta_N^2}, 
	& i=2,\cdots,N-1,\;k=0,\cdots,i-2, \\
\frac{4\theta_i\theta_{i-1}}{\theta_N^2}
	+ \frac{2\theta_{i-1}^2}{\theta_N^2}, 
	& i=1,\cdots,N-1,\;k=i-1, \\
\frac{2\theta_k}{\theta_N}, 
	& i=N,\;k=0,\cdots,i-2, \\
\frac{2\theta_{N-1}}{\theta_N} + \frac{2\theta_{N-1}^2}{\theta_N^2}, 
	& i=N,\;k=i-1, 
	\end{cases} \label{eq:sol0} \\
\hat{\lambda}_i &= 
	\frac{2\theta_{i-1}^2}{\theta_N^2}, \quad i=1,\cdots,N,
	\label{eq:sol1}  \\
\hat{\tau}_i &= \begin{cases}
	\frac{2\theta_i}{\theta_N^2},
	& i = 0,\cdots,N-1, \\
	1 - \frac{2\theta_{N-1}^2}{\theta_N^2} = \frac{1}{\theta_N},
	& i = N,
	\end{cases} \label{eq:sol2} \\
\hat{\gamma} &= 
	\frac{1}{\theta_N^2},
        \label{eq:sol3}
\end{align}
for
\begin{align}
\theta_i = \begin{cases}
        1, & i = 0, \\
        \frac{1 + \sqrt{1+4\theta_{i-1}^2}}{2}, & i = 1,\cdots,N-1, \\
	\frac{1 + \sqrt{1+8\theta_{i-1}^2}}{2}
	& i = N.
\end{cases}
\label{eq:thetai}
\end{align}
\begin{proof}
The following set of conditions
are sufficient for the feasible conditions of~\eqref{eq:HDR1}:
\begin{align}
\begin{cases}
\Qr(\bmr,\bmlam,\bmtau) 
	= 2\qr(\bmr,\bmlam,\bmtau)\qr(\bmr,\bmlam,\bmtau)^\top, & \\
\ttau = 2\qr(\bmr,\bmlam,\bmtau)\tau_N, & \\
\gamma = \tau_N^2, & \\
(\bmlam,\bmtau) \in \Lambda. & 
\end{cases}
\label{eq:feas}
\end{align}
The ``Appendix'' shows that the point $(\Hr,\bmhlam,\bmhtau,\hgam)$
in~\eqref{eq:sol0},~\eqref{eq:sol1},~\eqref{eq:sol2} and~\eqref{eq:sol3}
is the unique solution of~\eqref{eq:feas}
and also satisfies the feasible conditions of~\eqref{eq:HDR}.
\qed
\end{proof}
\end{lemma}

Note that the parameter $\theta_i$~\eqref{eq:thetai} 
used in Lemma~\ref{lem:feas}
differs from $t_i$~\eqref{eq:ti}
only at the last iteration $N$.
In other words, $\{\theta_0,\cdots,\theta_{N-1}\}$ is
equivalent to $\{t_0,\cdots,t_{N-1}\}$ in~\eqref{eq:ti}
satisfying~\eqref{eq:ti_sum},
whereas the last parameter $\theta_N$ satisfies
\begin{align}
\theta_N^2 - \theta_N - 2\theta_{N-1}^2 = 0
\label{eq:thetaN}
.\end{align}

The next lemma shows that the feasible point 
derived in Lemma~\ref{lem:feas} is
an optimal solution of both~\eqref{eq:HDR} and~\eqref{eq:HDR1}.
\begin{lemma}
\label{lem:solkkt}
The choice of $(\Hr,\bmhlam,\bmhtau,\hgam)$
in~\eqref{eq:sol0},~\eqref{eq:sol1},~\eqref{eq:sol2}
and~\eqref{eq:sol3}
is an optimal solution of both~\eqref{eq:HDR} and~\eqref{eq:HDR1}.
\begin{proof}
The proof in~\cite{kim:14:ofo-arxiv} uses 
the Karush-Kuhn-Tucker (KKT) conditions of linear
SDP~\eqref{eq:HDR}.
\qed
\end{proof}
\end{lemma}

The optimized step coefficients \Hh of interest are then derived 
using~\eqref{eq:hki,rki} 
with the analytical optimal solution $(\Hr,\bmhlam,\bmhtau,\hgam)$ of~\eqref{eq:HDR}.
It is interesting to note that the corresponding coefficients \Hh in~\eqref{eq:ogm2,h} below
have a recursive form that is similar to~\eqref{eq:fgm2,h} of FGM2,
as discussed further in Section~\ref{sec:opt,fo,form}.
\begin{lemma}
\label{lem:hd}
The choice of $(\Hh,\bmhlam,\bmhtau,\hgam)$ 
in~\eqref{eq:sol1},~\eqref{eq:sol2},~\eqref{eq:sol3} and
\begin{align}
\Hhkip
&= \begin{cases}
        \frac{1}{\theta_{i+1}}\paren{2\theta_k - \sum_{j=k+1}^i \Hhkj},
        & k = 0,\cdots,i-1, \\ 
        1 + \frac{2\theta_i - 1}{\theta_{i+1}},
        & k = i, 
\end{cases} \label{eq:ogm2,h}
\end{align}
for $i=0,\cdots,N-1$ with $\theta_i$ in~\eqref{eq:thetai}
is an optimal solution of~\eqref{eq:HD}.
\begin{proof}
Inserting \Hr~\eqref{eq:sol0}, \bmhlam~\eqref{eq:sol1} 
and \bmhtau~\eqref{eq:sol2} into~\eqref{eq:hki,rki},
and noting that $\hlam_i+\htau_i > 0$ for $i=1,\cdots,N$, we get
\begin{align*}
\Hhki 
&= \frac{\Hrki - \htau_i\sum_{j=k+1}^{i-1}\Hhkj}{\hlam_i + \htau_i},
	\quad i=1,\cdots,N,\;k=0,\cdots,i-1, \\
&= \begin{cases}
	\frac{\theta_N^2}{2\theta_i^2}\paren{
	\frac{4\theta_i\theta_k}{\theta_N^2} 
	- \frac{2\theta_i}{\theta_N^2}\sum_{j=k+1}^{i-1}\Hhkj}, 
	& i=1,\cdots,N-1,\; k=0,\cdots,i-2, \\ 
	\frac{\theta_N^2}{2\theta_i^2}\paren{
	\frac{4\theta_i\theta_{i-1}}{\theta_N^2}
	+ \frac{2\theta_{i-1}^2}{\theta_N^2}}
	= \frac{2\theta_i\theta_{i-1} + \theta_i^2 - \theta_i}{\theta_i^2}, 
	& i=1,\cdots,N-1,\; k=i-1, \\
	\frac{2\theta_k}{\theta_N} 
	- \frac{1}{\theta_N}\sum_{j=k+1}^{N-1}\Hhkj, 
	& i=N,\;k=0,\cdots,i-2, \\
	\frac{2\theta_{N-1}}{\theta_N} + \frac{2\theta_{N-1}^2}{\theta_N^2}
	= \frac{2\theta_N\theta_{N-1} + \theta_N^2 - \theta_N}{\theta_N^2}, 
	& i=N,\;k=i-1, 
\end{cases}
\end{align*}
which is equivalent to~\eqref{eq:ogm2,h}.
From~\cite[Theorem 3]{drori:14:pof},
the corresponding $(\Hh,\bmhlam,\bmhtau,\hgam)$
becomes an optimal solution of~\eqref{eq:HD}.
\qed
\end{proof}
\end{lemma}

The following theorem shows
that Algorithm FO with the optimized \Hh~\eqref{eq:ogm2,h}
achieves a new convergence bound.

\begin{theorem}
\label{thm:conv}
Let $f\;:\;\Reals^d\rightarrow\Reals$ be convex and $C_L^{1,1}$
and let $\bmx_0,\cdots,\bmx_N \in \Reals^d$ be generated by
Algorithm FO with \Hh~\eqref{eq:ogm2,h} for a given $N\ge1$. Then
\begin{align}
f(\bmx_N) - f(\bmx_*) \le \frac{L||\bmx_0 - \bmx_*||^2}{2\theta_N^2}
\le \frac{L||\bmx_0 - \bmx_*||^2}{(N+1)(N+1+\sqrt{2})},
\quad \forall \bmx_* \in X_*(f)
\label{eq:opt,fo,conv}
.\end{align}
\begin{proof}
Using \hgam~\eqref{eq:sol3} 
and $\theta_{N-1}^2 = t_{N-1}^2 \ge \frac{(N+1)^2}{4}$ 
from~\eqref{eq:ti_sum} and~\eqref{eq:thetai},
we get
\begin{align*}
\hat{\gamma}
& = \frac{1}{\theta_N^2}
	= \frac{4}{\paren{1 + \sqrt{1 + 8\theta_{N-1}^2}}^2}
	\le \frac{4}{\paren{1 + \sqrt{1 + 2(N+1)^2}}^2} \\
	& \le \frac{2}{(N+1)^2 + \sqrt{2}(N+1) + 1}
        \le \frac{2}{(N+1)(N+1+\sqrt{2})}
.\end{align*}
Then, we have
\begin{align*}
f(\bmx_N) - f(\bmx_*) \le  
\mathcal{B}_{\mathrm{D}}(\Hh,N,L,R)
= \frac{1}{2}LR^2\hat{\gamma}
\le \frac{LR^2}{(N+1)(N+1+\sqrt{2})},
\quad \forall \bmx_* \in X_*(f)
,\end{align*}
based on Lemma~\ref{lem:hd}.
Finally, we let $R = ||\bmx_0 - \bmx_*||$.
\qed
\end{proof}
\end{theorem}

Theorem~\ref{thm:conv} shows that algorithm FO
with the optimized \Hh~\eqref{eq:ogm2,h}
decreases the function $f$ with a bound that is twice
as small as that of Nesterov's fast gradient methods
in~\eqref{eq:fgm,conv} and~\eqref{eq:fgm,conv,aux},
confirming \DT's numerical results
in~\cite[Tables 1 and 2]{drori:14:pof}.
The proposed algorithm requires at most 
$N = \ceil{\sqrt{\frac{L}{\epsilon}}||\bmx_0 - \bmx_*||}$ iterations
to achieve the desired accuracy $f(\bmx_N) - f(\bmx_*) \le \epsilon$,
while Nesterov's fast gradient methods require
at most
$N = \ceil{\sqrt{\frac{2L}{\epsilon}}||\bmx_0 - \bmx_*||}$,
a factor of about $\sqrt{2}$-times more iterations.

The next section describes efficient implementations of 
the corresponding Algorithm FO with \Hh~\eqref{eq:ogm2,h}.

\section{Efficient formulations of proposed optimized first-order algorithms}
\label{sec:opt,fo,form}

Even though the analytical expression for \Hh in~\eqref{eq:ogm2,h} 
that solves~\eqref{eq:HD}
does not require an expensive SDP method,
using \Hh in Algorithm FO would still be computationally undesirable.
Noticing the similarity between~\eqref{eq:fgm2,h} of FGM2 and~\eqref{eq:ogm2,h},
we can expect that Algorithm FO with~\eqref{eq:ogm2,h}
may have an equivalent efficient form as FGM2,
as described next. 
In addition, we find an equivalent form of~\eqref{eq:ogm2,h}
that is similar to~\eqref{eq:fgm1,h} of FGM1,
so that we can find a formulation that is similar to FGM1
by analogy with how Proposition~\ref{prop:fgm1,fgm2} shows the equivalence
between~\eqref{eq:fgm1,h} and~\eqref{eq:fgm2,h}.

\begin{proposition}
The optimized \Hh in~\eqref{eq:ogm2,h} 
satisfies the following recursive relationship
\begin{align}
\Hhkip
&= \begin{cases}
	\frac{\theta_i - 1}{\theta_{i+1}} \Hhki, 
	& k = 0,\cdots,i-2, \\ 
	\frac{\theta_i - 1}{\theta_{i+1}}(\Hhimi - 1),
	& k = i - 1, \\
	1 + \frac{2\theta_i - 1}{\theta_{i+1}},
	& k = i,
\end{cases}
\label{eq:ogm1,h} 
\end{align}
for $i=0,\cdots,N-1$
with $\theta_i$ in~\eqref{eq:thetai}.
\begin{proof}
We follow the induction proof of Proposition~\ref{prop:fgm1,fgm2} 
showing the equivalence between~\eqref{eq:fgm1,h} and~\eqref{eq:fgm2,h}.
We use the notation $\Hhki'$ for the coefficient~\eqref{eq:ogm2,h}
to distinguish from~\eqref{eq:ogm1,h}.

It is obvious that
$\Hhiip' = \Hhiip,\; i=0,\cdots,N-1$,
and we clearly have
\begin{align*}
\Hhimip' &= \frac{1}{\theta_{i+1}}\paren{2\theta_{i-1} - \Hhimi'}
        = \frac{1}{\theta_{i+1}}\paren{2\theta_{i-1}
                - \paren{1 + \frac{2\theta_{i-1} - 1}{\theta_i}}} \\
        &= \frac{(2\theta_{i-1} - 1)(\theta_i - 1)}{\theta_i\theta_{i+1}}
        = \frac{\theta_i - 1}{\theta_{i+1}}\paren{\Hhimi - 1}
        = \Hhimip.
\end{align*}
for $i=0,\cdots,N-1$.

We next use induction by assuming $\Hhkip' = \Hhkip$
for $i=0,\cdots,n-1,\;k=0,\cdots,i$. We then have
\begin{align*}
\Hhknp' &= \frac{1}{\theta_{n+1}}
	\paren{2\theta_k - \sum_{j=k+1}^n \Hhkj'}
        = \frac{1}{\theta_{n+1}}\paren{2\theta_k - \sum_{j=k+1}^{n-1} \Hhkj'
                - \Hhkn'} \\
        &= \frac{\theta_n - 1}{\theta_{n+1}}\Hhkn'
	= \frac{\theta_n - 1}{\theta_{n+1}}\Hhkn 
	= \Hhknp
\end{align*}
for $k=1,\cdots,n-2$.
Note that this proof is independent of the choice of $\theta_i$.
\qed
\end{proof}
\end{proposition}

Next, we revisit the derivation in Section~\ref{sec:fo,ex}
to transform Algorithm FO with~\eqref{eq:ogm2,h} or~\eqref{eq:ogm1,h}
into efficient formulations
akin to Nesterov's fast gradient methods,
leading to practical algorithms.

\subsection{Proposed optimized gradient method 1 (OGM1)}

We first propose the following optimized gradient method, called OGM1, 
using~\eqref{eq:ogm1,h} in Algorithm FO.
OGM1 is computationally similar to FGM1
yet the sequence $\{\bmx_i\}$ generated by OGM1
achieves the fast convergence bound in Theorem~\ref{thm:conv}.

\fbox{
\begin{minipage}[t]{0.85\linewidth}
\vspace{-10pt}
\begin{flalign*}
&\quad \text{\bf Algorithm OGM1} & \\
&\qquad \text{Input: } f\in C_L^{1,1}(\Reals^d)\text{ convex},\; \bmx_0\in\Reals^d,\;
	\bmy_0 = \bmx_0,\; \theta_0 = 1. & \\
&\qquad \text{For } i = 0,\cdots,N-1 & \\
&\qquad \qquad \bmy_{i+1} = \bmx_i - \frac{1}{L}f'(\bmx_i) & \\
&\qquad \qquad \theta_{i+1}
                = \begin{cases}
			\frac{1+\sqrt{1+4\theta_i^2}}{2}, & i \le N-2 \\
			\frac{1 + \sqrt{1+8\theta_i^2}}{2}, & i = N-1
			\end{cases} & \\
&\qquad \qquad \bmx_{i+1} = \bmy_{i+1}
                + \frac{\theta_i - 1}{\theta_{i+1}}(\bmy_{i+1} - \bmy_i) 
		+ \frac{\theta_i}{\theta_{i+1}}(\bmy_{i+1} - \bmx_i)& 
\end{flalign*}
\end{minipage}
} \vspace{5pt}

\noindent
Apparently, the proposed OGM1 accelerates FGM1 by
using just one additional momentum term 
$\frac{\theta_i}{\theta_{i+1}}(\bmy_{i+1} - \bmx_i)$,
and thus OGM1 is computationally efficient.
Also, unlike \DT's approach that requires choosing $N$ 
for using a SDP solver before iterating,
the proposed OGM1 does not need to know $N$ in advance
because the coefficients \Hh (or $\theta_i$) for 
intermediate iterations ($i=0,\cdots,N-1$) do not depend on $N$.

\begin{proposition}
The sequence $\{\bmx_0,\cdots,\bmx_N\}$ generated by Algorithm FO
with~\eqref{eq:ogm1,h} is identical
to the corresponding sequence generated by Algorithm OGM1.
\begin{proof}
We use induction,
and for clarity, we use the notation $\bmx_0',\cdots,\bmx_N'$
for Algorithm FO.
It is obvious that $\bmx_0' = \bmx_0$, and since $\theta_0 = 1$ we get
\begin{align*}
\bmx_1' = \bmx_0' - \frac{1}{L}\hat{h}_{1,0}f'(\bmx_0')
	= \bmx_0 - \frac{1}{L}\paren{1 + \frac{2\theta_0 - 1}{\theta_1}}f'(\bmx_0)
	= \bmy_1 + \frac{\theta_0}{\theta_1}(\bmy_1 - \bmx_0)
	= \bmx_1
.\end{align*}
Assuming $\bmx_i' = \bmx_i$ for $i=0,\cdots,n$, we then have
{\allowdisplaybreaks[4]
\begin{align*}
\bmx_{n+1}' 
	=& \bmx_n' - \frac{1}{L}\Hhnnp f'(\bmx_n') 
	- \frac{1}{L}\Hhnmnp f'(\bmx_{n-1}')
	- \frac{1}{L}\sum_{k=0}^{n-2}\Hhknp f'(\bmx_k') \\
	=& \bmx_n - \frac{1}{L}\paren{1 + \frac{2\theta_n - 1}{\theta_{n+1}}}f'(\bmx_n) \\ 
	&- \frac{\theta_n - 1}{\theta_{n+1}}(\Hhnmn - 1)f'(\bmx_{n-1}) 
	- \frac{1}{L}\sum_{k=0}^{n-2}
	\frac{\theta_n - 1}{\theta_{n+1}}\Hhkn f'(\bmx_k) \\
	=& \bmx_n - \frac{1}{L}\paren{1 + \frac{\theta_n}{\theta_{n+1}}}f'(\bmx_n) \\
	&+ \frac{\theta_n-1}{\theta_{n+1}}\paren{-\frac{1}{L}f'(\bmx_n)
	+ \frac{1}{L}f'(\bmx_{n-1}) 
	- \frac{1}{L}\sum_{k=0}^{n-1}\Hhkn f'(\bmx_k)} \\
	=& \bmy_{n+1} + \frac{\theta_n}{\theta_{n+1}}(\bmy_{n+1} - \bmx_n) 
	+ \frac{\theta_n-1}{\theta_{n+1}}\paren{-\frac{1}{L}f'(\bmx_n)
	+ \frac{1}{L}f'(\bmx_{n-1})
	+ \bmx_n - \bmx_{n-1}} \\
	=& \bmy_{n+1} + \frac{\theta_n - 1}{\theta_{n+1}}(\bmy_{n+1} - \bmy_n)
		+ \frac{\theta_n}{\theta_{n+1}}(\bmy_{n+1} - \bmx_n)
	= \bmx_{n+1}.
\end{align*}}
\qed
\end{proof}
\end{proposition}

\subsection{Proposed optimized gradient method 2 (OGM2)}

We propose another efficient formulation
of Algorithm FO with~\eqref{eq:ogm2,h}
that is similar to the formulation of FGM2.

\fbox{
\begin{minipage}[t]{0.85\linewidth}
\vspace{-10pt}
\begin{flalign*}
&\quad \text{\bf Algorithm OGM2} & \\
&\qquad \text{Input: } f\in C_L^{1,1}(\Reals^d)\text{ convex},\; \bmx_0\in\Reals^d,\;
	 \theta_0 = 1. & \\
&\qquad \text{For } i = 0,\cdots,N-1 & \\
&\qquad \qquad \bmy_{i+1} = \bmx_i - \frac{1}{L}f'(\bmx_i) & \\
&\qquad \qquad \bmz_{i+1} = \bmx_0 - \frac{1}{L}\sum_{k=0}^i 2\theta_k f'(\bmx_k) & \\
&\qquad \qquad \theta_{i+1} = \begin{cases}
			\frac{1+\sqrt{1+4\theta_i^2}}{2}, & i \le N-2 \\
			\frac{1 + \sqrt{1+8\theta_{i}^2}}{2}, & i = N-1 
			\end{cases} & \\
&\qquad \qquad \bmx_{i+1} = \paren{1 - \frac{1}{\theta_{i+1}}}\bmy_{i+1} 
				+ \frac{1}{\theta_{i+1}}\bmz_{i+1} & 
\end{flalign*}
\end{minipage}
} \vspace{5pt}

\noindent
The sequence $\{\bmx_i\}$ generated by OGM2 achieves 
the fast convergence bound in Theorem~\ref{thm:conv}.
Algorithm OGM2 doubles the weight on 
all previous gradients for $\{\bmz_i\}$ compared to FGM2,
providing some intuition for its two-fold acceleration.
OGM2 requires comparable computation per iteration as FGM2.

\begin{proposition}
The sequence $\{\bmx_0,\cdots,\bmx_N\}$ generated by Algorithm FO
with~\eqref{eq:ogm2,h} is identical
to the corresponding sequence generated by Algorithm OGM2.
\begin{proof}
We use induction,
and for clarity, we use the notation $\bmx_0',\cdots,\bmx_N'$
for Algorithm FO.
It is obvious that $\bmx_0' = \bmx_0$, and since $\theta_0 = 1$ we get
\begin{align*}
\bmx_1' = \bmx_0' - \frac{1}{L}\hat{h}_{1,0}f'(\bmx_0')
        = \bmx_0 - \frac{1}{L}\paren{1 + \frac{2\theta_0 - 1}{\theta_1}}f'(\bmx_0)
        = \bmy_1 + \frac{\theta_0}{\theta_1}(\bmy_1 - \bmx_0)
        = \bmx_1
.\end{align*}
Assuming $\bmx_i' = \bmx_i$ for $i=0,\cdots,n$, we then have
{\allowdisplaybreaks[4]
\begin{align*}
\bmx_{n+1}' 
=& \bmx_n' - \frac{1}{L}\Hhnnp f'(\bmx_n') 
	- \frac{1}{L}\sum_{k=0}^{n-1} \Hhknp f'(\bmx_k') \\
=& \bmx_n - \frac{1}{L}\paren{1 + \frac{2\theta_n - 1}{\theta_{n+1}}} f'(\bmx_n)
	- \frac{1}{L}\sum_{k=0}^{n-1} \frac{1}{\theta_{n+1}}
		\paren{2\theta_k - \sum_{j=k+1}^n \Hhkj}f'(\bmx_k) \\
=& \paren{1 - \frac{1}{\theta_{n+1}}}\paren{\bmx_n - \frac{1}{L}f'(\bmx_n)} 
	+ \frac{1}{\theta_{n+1}}\paren{\bmx_0 - \frac{1}{L}\sum_{k=0}^n 2\theta_kf'(\bmx_k)} \\
=& \paren{1 - \frac{1}{\theta_{n+1}}}\bmy_{n+1}
	+ \frac{1}{\theta_{n+1}}\bmz_{n+1}
= \bmx_{n+1}.
\end{align*}}
$\!\!\!$The third  
equality uses the telescoping sum
$\bmx_n = \bmx_0 + \sum_{j=1}^{n} (\bmx_j - \bmx_{j-1})$
and~\eqref{eq:fo} in Algorithm FO.
\qed
\end{proof}
\end{proposition}

\section{Discussion}
\label{sec:disc}

After submitting this work~\cite{kim:14:ofo-arxiv}, 
Taylor \etal~\cite{taylor:15:ssc} further studied the PEP approach
to compute the exact worst-case bound of first-order methods,
unlike \DT~\cite{drori:14:pof} and this paper that use the \emph{relaxed} PEP.
Taylor \etal~\cite{taylor:15:ssc}
studied the tightness of relaxations on PEP introduced in~\cite{drori:14:pof}
and avoided some strict relaxations.

Inspired by~\cite[Conjecture 5]{taylor:15:ssc},
we developed the following theorem
that shows that the smallest upper bound in~\eqref{eq:opt,fo,conv} for OGM1 and OGM2 is tight,
despite the various relaxations of PEP
used in~\cite{drori:14:pof} and herein.
(Similar tightness results are shown for the gradient methods
with a constant step size $0 < h \le 1$ in~\cite{drori:14:pof}.)
The following theorem specifies a worst-case convex function $\phi(\bmx)$
in $C_L^{1,1}(\Reals^d)$
for which the optimized gradient methods achieve 
their smallest upper bound in~\eqref{eq:opt,fo,conv}.

\begin{theorem}
\label{thm:low}
For the following convex functions in $C_L^{1,1}(\Reals^d)$
for all $d\ge1$: 
\begin{align}
\phi(\bmx) = \begin{cases}
        \frac{LR}{\theta_N^2}||\bmx|| - \frac{LR^2}{2\theta_N^4}, 
                & \text{if } ||\bmx||\ge\frac{R}{\theta_N^2}, \\
        \frac{L}{2}||\bmx||^2, 
                & \text{otherwise},
\end{cases}
\label{eq:worst}
\end{align}
both OGM1 and OGM2 exactly achieve 
the smallest upper bound in~\eqref{eq:opt,fo,conv}, \ie,
\begin{align*}
\phi(\bmx_N) - \phi(\bmx_*) = \frac{L||\bmx_0 - \bmx_*||^2}{2\theta_N^2}
.\end{align*}
\begin{proof}
We show in the ``Appendix'' that the following property 
of the coefficients \Hh~\eqref{eq:ogm1,h} of OGM1 and OGM2 holds:
\begin{align}
\sum_{j=1}^i\sum_{k=0}^{j-1}\Hhkj = \begin{cases}
        \theta_i^2 - 1, & i=1,\cdots,N-1, \\
        \frac{1}{2}(\theta_N^2 - 1), & i = N.
        \end{cases}
\label{eq:ogm,h,prop}
\end{align}

Then, starting from $\bmx_0 = R\bmnu$, where \bmnu is a unit vector,
and using \eqref{eq:ogm,h,prop},
the iterates of OGM1 and OGM2 are as follows
\begin{align*}
\bmx_i &= \bmx_0 - \frac{1}{L}\sum_{j=1}^i\sum_{k=0}^{j-1}\Hhkj \phi'(\bmx_k) =
	\begin{cases}
		\paren{1 - \frac{\theta_i^2 - 1}{\theta_N^2}}R\bmnu, & i=0,\cdots,N-1, \\
		\paren{1 - \frac{\theta_N^2 - 1}{2\theta_N^2}}R\bmnu, & i=N, \\		
	\end{cases}
\end{align*}
where the corresponding sequence $\{\bmx_0,\cdots,\bmx_N\}$ stays 
in the affine region of the function $\phi(\bmx)$
with the same gradient value:
\begin{align*}
\phi'(\bmx_i) &= \frac{LR}{\theta_N^2}\bmnu, \quad i=0,\cdots,N.
\end{align*}
Therefore, after $N$ iterations of OGM1 and OGM2, we have
\begin{align*}
\phi(\bmx_N) - \phi(\bmx_*) = \phi(\bmx_N) = \frac{LR^2}{2\theta_N^2} 
,\end{align*}
exactly matching the smallest upper bound in~\eqref{eq:opt,fo,conv}.
\qed
\end{proof}
\end{theorem}

This result implies that the exact PEP bound $\mathcal{B}_{\mathrm{P}}(\Hh,N,\dd,L,R)$ 
of OGM1 and OGM2 
is equivalent to their relaxed bound $\mathcal{B}_{\mathrm{D}}(\Hh,N,L,R)$
that is independent of $d$.
Note that Taylor \etal~\cite{taylor:15:ssc} showed 
that the exact PEP $\mathcal{B}_{\mathrm{P}}(\bmh,N,\dd,L,R)$
is independent of $d$ 
{\color{blue}
for the large-scale condition ``$d \ge N + 2$''.
}
Whereas the OGM bound~\eqref{eq:opt,fo,conv} is tight,
the FGM bounds~\eqref{eq:fgm,conv} and~\eqref{eq:fgm,conv,aux}
are not tight~\cite{drori:14:pof,taylor:15:ssc},
somewhat weakening the utility of the fact
that the OGM bound~\eqref{eq:opt,fo,conv} is twice smaller than the FGM bounds.
However, 
Figure 5 in~\cite{taylor:15:ssc} shows
that the FGM bounds~\eqref{eq:fgm,conv} and~\eqref{eq:fgm,conv,aux} 
become close to tight asymptotically
as $N$ increases,
so the factor of $2$ can have practical value
when using many iterations.
We leave more complete comparisons as future work.

\section{Conclusion}
\label{sec:conc}

We proposed new optimized first-order algorithms
that achieve a worst-case convergence bound that is
twice as small as that of Nesterov's methods
for smooth unconstrained convex minimization,
inspired by Drori and Teboulle~\cite{drori:14:pof}.
The proposed first-order methods 
are comparably efficient for implementation as Nesterov's methods.
Thus it is natural to use the proposed OGM1 and OGM2 to replace Nesterov's methods 
in smooth unconstrained convex minimization. 
Numerical results in large-scale imaging applications show
practical convergence acceleration
consistent with those predicted by the bounds given here
\cite{kim:14:oms,kim:15:aof}.
Those applications use regularizers
that have shapes somewhat similar 
to the worst-case function~\eqref{eq:worst}.

The efficient formulations of both Nesterov's methods and the new optimized first-order methods
still seem somewhat magical.
Recently, \cite{allen-zhu:15:lca}, \cite{odonoghue:15:arf} and \cite{su:15:ade}
studied Nesterov's FGM formulations,
and extending such studies to the new OGM methods should further illuminate
the fundamental causes for
their efficient formulations and acceleration.
Also, new optimized first-order methods lack analytical convergence bounds
for the intermediate iterations,
whereas numerical bounds are studied in~\cite{taylor:15:ssc};
deriving those analytical bounds is interesting future work.

Drori recently extended the PEP approach to projected gradient methods
for constrained smooth convex minimization~\cite{drori:14:phd}.
Extending this approach to general first-order algorithms including our proposed OGM1 and OGM2
is important future work.
In addition, just as Nesterov's fast gradient methods
have been extended for nonsmooth composite convex minimization
\cite{beck:09:afi,nesterov:13:gmf},
extending the proposed optimized first-order algorithms
for minimizing nonsmooth composite convex functions
would be a natural direction to pursue.

While \DT's PEP approach involves a series of relaxations
to make the problem solvable,
OGM1 and OGM2 with the step coefficients \Hh
that are optimized over the \emph{relaxed} PEP upper bound~\eqref{eq:HD} 
achieve an \emph{exact} bound in Theorem~\ref{thm:low}.
However,
it remains an open problem to either prove  
that the smallest upper bound in~\eqref{eq:opt,fo,conv} of OGM1 and OGM2 is optimal.
We leave either proving the above statement for~\eqref{eq:HP}
or to further optimize the first-order methods
as future work. 

\section{Appendix}
\subsection{Proof of Lemma~\ref{lem:feas}}

We prove that the choice $(\Hr,\bmhlam,\bmhtau,\gam)$
in~\eqref{eq:sol0},~\eqref{eq:sol1},~\eqref{eq:sol2} and~\eqref{eq:sol3}
satisfies the feasible conditions~\eqref{eq:feas}
of~\eqref{eq:HDR1}.

Using the definition of
$\Qr(\bmr,\bmlam,\bmtau)$ in~\eqref{eq:Qr},
and considering the first two conditions of~\eqref{eq:feas}, we get
\begin{align*}
\lambda_{i+1} &=
\breve{Q}_{i,i}(\bmr,\bmlam,\bmtau)
        = 2\breve{q}_i^2(\bmr,\bmlam,\bmtau) 
	= \frac{1}{2\tau_N^2}\tau_i^2 \\
	&= \begin{cases}
                        \frac{1}{2(1-\lambda_N)^2}\lambda_1^2, & i = 0 \\
                        \frac{1}{2(1-\lambda_N)^2}
                                (\lambda_{i+1} - \lambda_i)^2,
                        & i = 1,\cdots,N-1,
                \end{cases}
\end{align*}
where the last equality comes from
$(\bmlam,\bmtau)\in\Lambda$,
and this reduces to the following recursion:
\begin{align}
\begin{cases}
\lambda_1 = 2(1 - \lambda_N)^2, & \\
(\lambda_i - \lambda_{i-1})^2 - \lambda_1\lambda_i = 0. \quad i=2,\cdots,N. &
\end{cases}
\label{eq:recur}
\end{align}

We use induction to prove that
the solution of~\eqref{eq:recur} is
\begin{align*}
\lambda_i =
\begin{cases}
\frac{2}{\theta_N^2}, & i=1, \\
\theta_{i-1}^2\lambda_1, & i=2,\cdots,N,
\end{cases}
\end{align*}
which is equivalent to \bmhlam~\eqref{eq:sol1}.
It is obvious that $\lambda_1 = \theta_0\lambda_1$,
and for $i=2$ in~\eqref{eq:recur}, we get
\begin{align*}
\lambda_2 = \frac{3\lambda_1 + \sqrt{9\lambda_1^2 - 4\lambda_1^2}}{2}
        = \frac{3+\sqrt{5}}{2}\lambda_1 = \theta_1^2\lambda_1
.\end{align*}
Then, assuming $\lambda_i = \theta_{i-1}^2\lambda_1$ for $i=1,\cdots,n$
and $n\le N-1$,
and using the second equality in~\eqref{eq:recur} for $i=n+1$, we get
\begin{align*}
\lambda_{n+1} &= \frac{\lambda_1 + 2\lambda_n
        + \sqrt{(\lambda_1 + 2\lambda_n)^2 - 4\lambda_n^2}}{2}
        = \frac{1 + 2\theta_{n-1}^2 + \sqrt{1 + 4\theta_{n-1}^2}}{2}\lambda_1 \\
        &= \paren{\theta_{n-1}^2
                + \frac{1 + \sqrt{1 + 4\theta_{n-1}^2}}{2}}\lambda_1
        = \theta_n^2\lambda_1,
\end{align*}
where the last equality uses~\eqref{eq:ti_sum}.
Then we use the first equality in~\eqref{eq:recur}
to find the value of $\lambda_1$ as
\begin{align*}
& \lambda_1 = 2(1 - \theta_{N-1}^2\lambda_1)^2 \\
& \theta_{N-1}^4\lambda_1^2 - 2\paren{\theta_{N-1}^2 + \frac{1}{4}}\lambda_1 + 1 = 0 \\
& \lambda_1
        = \frac{\theta_{N-1}^2 + \frac{1}{4}
                - \sqrt{(\theta_{N-1}^2 + \frac{1}{4})^2 - \theta_{N-1}^4}}{\theta_{N-1}^4}
        = \frac{1}{\theta_{N-1}^2 + \frac{1}{4} + \sqrt{\frac{\theta_{N-1}^2}{2} + \frac{1}{16}}} \\
&\quad = \frac{8}{\paren{1 + \sqrt{1 + 8\theta_{N-1}^2}}^2}
        = \frac{2}{\theta_N^2}
\end{align*}
with $\theta_N$ in~\eqref{eq:thetai}.

Until now, we derived $\bmhlam$~\eqref{eq:sol1}
using some conditions of~\eqref{eq:feas}.
Consequently, using the last two conditions in~\eqref{eq:feas}
with~\eqref{eq:ti_sum} and~\eqref{eq:thetaN},
we can easily derive the following:
\begin{align*}
\tau_i &= \begin{cases}
\hlam_1 = \frac{2}{\theta_N^2}, & i = 0, \\
\hlam_{i+1} - \hlam_i = \frac{2\theta_i^2}{\theta_N^2}
                - \frac{2\theta_{i-1}^2}{\theta_N^2}
                = \frac{2\theta_i}{\theta_N^2}, & i=1,\cdots,N-1, \\
1 - \hlam_N = 1 - \frac{2\theta_{N-1}^2}{\theta_N^2}
        = \frac{1}{\theta_N}, & i = N,
\end{cases} \\
\gamma &= \tau_N^2 = \frac{1}{\theta_N^2},
\end{align*}
which are equivalent to \bmhtau~\eqref{eq:sol2} and \hgam~\eqref{eq:sol3}.

Next, we derive \Hr for given 
\bmhlam~\eqref{eq:sol1} and \bmhtau~\eqref{eq:sol2}.
Inserting \bmhtau~\eqref{eq:sol2} to 
the first two conditions of~\eqref{eq:feas},
we get
\begin{align}
\begin{cases}
\breve{q}_i(\Hr,\bmhlam,\bmhtau) 
	= \frac{\hat{\tau}_i}{2\hat{\tau}_N}
	= \frac{\theta_i}{\theta_N}, &  \\
\breve{Q}_{i,k}(\Hr,\bmhlam,\bmhtau)
        = 2\breve{q}_i(\bmr,\bmhlam,\bmhtau)
                \breve{q}_k(\bmr,\bmhlam,\bmhtau)
        = \frac{2\theta_i\theta_k}{\theta_N^2}, &
\end{cases}
\label{eq:Qrqr}
\end{align}
for $i,k=0,\cdots,N-1$,
and considering~\eqref{eq:notation} and~\eqref{eq:Qrqr}, we get
\begin{align}
\breve{S}_{i,k}(\Hr,\bmhlam,\bmhtau)
&= \begin{cases}
        \frac{2\theta_i\theta_k}{\theta_N^2},        
		& i,k = 0,\cdots,N-1, \\ 
        \frac{\theta_i}{\theta_N}        
		& i=0,\cdots,N-1,\;k=N, \\ 
        \frac{\theta_k}{\theta_N},
                & i = N,\;k=0,\cdots,N-1, \\ 
        \frac{1}{2}
                & i = N,\;k = N.
        \end{cases}
\label{eq:S1}
\end{align}
Finally, using the two equivalent forms~\eqref{eq:Sr} and~\eqref{eq:S1}
of $\Sr(\Hr,\bmhlam,\bmhtau)$, we get
\begin{align}
\breve{S}_{i,k}(\Hr,\bmhlam,\bmhtau)
= 
\begin{cases}
\frac{1}{2}\Hrki 
	= \frac{2\theta_i\theta_k}{\theta_N^2},
	& i=2,\cdots,N-1,\;k=0,\cdots,i-2, \\
\frac{1}{2}(\Hrki - \hlam_i)
	= \frac{2\theta_i\theta_k}{\theta_N^2},
	& i=1,\cdots,N-1,\;k=i-1, \\
\frac{1}{2}\Hrki
	= \frac{\theta_k}{\theta_N},
	& i=N,\;k=0,\cdots,i-2, \\
\frac{1}{2}(\Hrki - \hlam_i)
	= \frac{\theta_k}{\theta_N}.
        & i=N,\;k=i-1, \\
\end{cases}
\end{align}
and this can be easily converted 
to the choice \Hrki in~\eqref{eq:sol0}.

For these given $(\Hr,\bmhlam,\bmhtau)$, we can easily notice that
\begin{align}
\paren{\begin{array}{cc}
        \Sr(\Hr,\bmhlam,\bmhtau) & \frac{1}{2}\bmhtau \\
        \frac{1}{2}\bmhtau^\top & \frac{1}{2}\hgam
\end{array}}
= \paren{\begin{array}{cc} \frac{2}{\theta_N^2}\bmttheta\bmttheta^\top
                & \frac{1}{\theta_N^2}\bmttheta \\
                \frac{1}{\theta_N^2}\bmttheta^\top
                & \frac{1}{2\theta_N^2} \end{array}} 
= \frac{2}{\theta_N^2}\paren{\begin{array}{c} \bmttheta \\ \frac{1}{2} \end{array}}
                \paren{\begin{array}{c} \bmttheta \\ \frac{1}{2} \end{array}}^\top
                \succeq 0
\label{eq:r,feas}
\end{align}
for $\bmttheta = \paren{\theta_0,\cdots,\theta_{N-1},\frac{\theta_N}{2}}^\top$,
showing that the choice is feasible in both~\eqref{eq:HDR} and~\eqref{eq:HDR1}.
\qed

\subsection{Proof of~\eqref{eq:ogm,h,prop}}

We prove that~\eqref{eq:ogm,h,prop} holds for
the coefficients \Hh~\eqref{eq:ogm1,h} of OGM1 and OGM2.

We first show the following property using induction:
\begin{align*}
\sum_{k=0}^{j-1}\Hhkj = \begin{cases}
	\theta_j, & j = 1,\cdots,N-1, \\
	\frac{1}{2}(\theta_N+1), & j = N.
\end{cases}
\end{align*}
Clearly, $\Hhzo = 1 + \frac{2\theta_0 - 1}{\theta_1} = \theta_1$ using~\eqref{eq:ti_sum}.
Assuming $\sum_{k=0}^{j-1}\Hhkj = \theta_j$ for $j=1,\cdots,n$ and $n\le N-1$, we get
\begin{align*}
\sum_{k=0}^n\Hhknp 
	&= 1 + \frac{2\theta_n - 1}{\theta_{n+1}}
		+ \frac{\theta_n - 1}{\theta_{n+1}}(\Hhnmn - 1)
		+ \frac{\theta_n - 1}{\theta_{n+1}}\sum_{k=0}^{n-2}\Hhkn \\
	&= 1 + \frac{\theta_n}{\theta_{n+1}} 
		+ \frac{\theta_n - 1}{\theta_{n+1}}\sum_{k=0}^{n-1}\Hhkn
	= \frac{\theta_{n+1} + \theta_n^2}{\theta_{n+1}} \\
	&= \begin{cases}
		\theta_n, & n = 1,\cdots,N-2, \\
		 \frac{1}{2}(\theta_N + 1), & n = N-1,	
	\end{cases} 
\end{align*}
where the last equality uses~\eqref{eq:ti_sum} and~\eqref{eq:thetaN}.

Then,~\eqref{eq:ogm,h,prop} can be easily derived
using~\eqref{eq:ti_sum} and~\eqref{eq:thetaN} as
\begin{align*}
\sum_{j=1}^i\sum_{k=0}^{j-1}\Hhkj 
	&= \begin{cases}
		\sum_{j=1}^i\theta_j, & i = 1,\cdots,N-1, \\
		\sum_{j=1}^{N-1}\theta_j + \frac{1}{2}(\theta_N+1), & i = N,
	\end{cases} \\
	&= \begin{cases}
        	\theta_i^2 - 1, & i = 1,\cdots,N-1, \\
        	\frac{1}{2}(\theta_N^2 - 1), & i = N.
	\end{cases}
\end{align*}
\qed

\bibliographystyle{spmpsci}      
\bibliography{master0,mastersub}   

%
%
%

\end{document}